\documentclass[journal]{IEEEtran}
\ifCLASSINFOpdf
\else
\fi

\usepackage{lipsum}
\usepackage{epstopdf}
\usepackage{amsfonts}
\usepackage{mathrsfs} 
\usepackage{amsmath} 
\usepackage{amssymb} % JA added
\usepackage{psfrag} 
\usepackage{enumerate}
\usepackage{mathtools}
\usepackage{algorithm}
\usepackage[noend]{algorithmic}
\usepackage{scalerel}
\usepackage{comment}
\usepackage{gensymb}
\usepackage{tabularx}
\usepackage{multirow}
\usepackage{subcaption}
\usepackage{accents}
\usepackage{setspace}
\usepackage{xcolor}

\ifpdf
  \DeclareGraphicsExtensions{.eps,.pdf,.png,.jpg}
\else
  \DeclareGraphicsExtensions{.eps}
\fi

\newtheorem{assumption}{Assumption}
\newcolumntype{Y}{>{\centering\arraybackslash}X}

\newcommand{\dbtilde}[1]{\accentset{\approx}{#1}}

\DeclareMathOperator*{\argmin}{arg\,min}
\DeclareMathOperator*{\mind}{\min.}
\newcommand{\minimize}{\mbox{minimize}}
\DeclareMathOperator{\sgn}{sgn}

\newcommand{\st}{\text{subject to}}
\newcommand{\suchthat}{\text{such that}}
\newcommand{\Real}{\mathbb{R}}

\newcommand{\Sspace}{\mathbb{S}^+}
\newcommand{\Zset}{\mathbb{Z}}
\newcommand{\T}{\mathsf{T}}

\newcommand{\rank}{{\rm rank}}
\newcommand{\inte}{{\mathrm{int}}}
\newcommand{\clos}{{\mathrm{clos}}}
\newcommand{\base}{{\bf{e}}}
\newcommand{\eye}{\mathbf{I}}
\newcommand{\zero}{\mathbf{0}}
\newcommand{\one}{\mathbf{1}}

\newcommand{\fea}{{\rm fea}}
\newcommand{\diff}{\mathrm{d}}

\newcommand{\alphavec}{\boldsymbol{\alpha}}
\newcommand{\betavec}{\boldsymbol{\beta}}

\newcommand{\Graph}{\mathcal{G}}

\newcommand{\Vertex}{\mathcal{V}}
\newcommand{\VertexG}{\mathcal{V}_{\rm{G}}}
\newcommand{\VertexL}{\mathcal{V}_{\rm{L}}}
\newcommand{\Edge}{\mathcal{E}}

\newcommand{\nG}{N_{\rm{G}}}
\newcommand{\nL}{N_{\rm{L}}}
\newcommand{\nGL}{N}
\newcommand{\nE}{E}

\newcommand{\CM}{\mathbf{C}}
\newcommand{\B}{\mathbf{B}}

\newcommand{\sgen}{\mathbf{s}^g}
\newcommand{\sgens}{\mathbf{s}_*^g}
\newcommand{\sload}{\mathbf{s}^l}
\newcommand{\sloads}{\mathbf{s}_*^l}

\newcommand{\pf}{\mathbf{p}}
\newcommand{\genllim}{\underline{\mathbf{s}}^g}
\newcommand{\genulim}{\overline{\mathbf{s}}^g}

\newcommand{\ang}{\boldsymbol{\theta}}
\newcommand{\angs}{\boldsymbol{\theta}_*}
\newcommand{\f}{{\bf f}}
\newcommand{\fs}{{\bf f}_*}
\newcommand{\fss}{{\bf f}_{**}}
\newcommand{\M}{{\bf M}}
\newcommand{\RelationMat}{{\bf H}}
\newcommand{\subRelationMat}{{\bf H}_1}
\newcommand{\bigMat}{{\bf R}}
\newcommand{\pflow}{{\bf p}}
\newcommand{\pflows}{{\bf p}_*}
\newcommand{\taueq}{\boldsymbol{\tau}}
\newcommand{\mup}{\boldsymbol{\mu}_+}
\newcommand{\mum}{\boldsymbol{\mu}_-}
\newcommand{\mupm}{\boldsymbol{\mu}}
\newcommand{\lambdap}{\boldsymbol{\lambda}_+}
\newcommand{\lambdam}{\boldsymbol{\lambda}_-}
\newcommand{\lambdapm}{\boldsymbol{\lambda}}
\newcommand{\set}{\Omega}
\newcommand{\setbad}{\Omega_{\rm bad}}
\newcommand{\setf}{\Omega_{\f}}
\newcommand{\setsl}{\Omega_{\sload}}

\newcommand{\setpara}{\Omega_{\para}}
\newcommand{\setall}{\Omega}
\newcommand{\setallr}{\tilde{\Omega}}
\newcommand{\para}{{\boldsymbol{\xi}}}
\newcommand{\paras}{{\boldsymbol{\xi}}_*}
\newcommand{\parass}{{\boldsymbol{\xi}}_{**}}
\newcommand{\setslr}{\widetilde{\Omega}_{\sload}}

\newcommand{\setparar}{\widetilde{\Omega}_{\para}}

\newcommand{\lbold}{\boldsymbol{\lambda}}
\newcommand{\nubold}{\boldsymbol{\nu}}

\newcommand{\setR}{\mathcal{R}}
\newcommand{\setS}{\mathcal{S}}
\newcommand{\setT}{\mathcal{T}}
\newcommand{\setQ}{\mathcal{Q}}

\newcommand{\OPF}{\mathcal{OPF}}

\newcommand{\JM}{\mathbf{J}}

\newcommand{\setSgen}{\setS_{\rm G}}
\newcommand{\setSbra}{\setS_{\rm B}}

\newcommand{\dD}{\mathrm{d}}

\newcommand{\x}{{\bf x}}
\newcommand{\y}{{\bf y}}
\newcommand{\X}{{\bf X}}
\newcommand{\Y}{{\bf Y}}

\newcommand{\z}{{\bf z}}
\newcommand{\Aeq}{{\bf A}_{\rm{eq}}}
\newcommand{\Ain}{{\bf A}_{\rm{in}}}
\newcommand{\beq}{{\bf b}_{\rm{eq}}}
\newcommand{\bin}{{\bf b}_{\rm{in}}}
\newcommand{\bpara}{{\bf b}'}
\newcommand{\Aintilde}{\tilde{{\bf A}}_{\rm{in}}}
\newcommand{\Aindbtilde}{\dbtilde{{\bf A}}_{\rm{in}}}
\newcommand{\setbinding}{\mathcal{X}}

\newcommand{\sgenh}{\mathbf{\hat{s}}^g}
\newcommand{\sloadh}{\mathbf{\hat{s}}^l}
\newcommand{\pflowh}{\mathbf{\hat{\pflow}}}
\newcommand{\angh}{\boldsymbol{\hat{\theta}}}

\newcommand{\dir}{\mathbf{v}}

\newcommand{\sen}{C}
\newcommand{\wcsen}{C^{\rm wc}}
\newcommand{\ifromj}{{i\leftarrow j}}

\newcommand{\sign}{\textbf{sign}}
\newcommand{\trace}{\textbf{tr}}

\newcommand{\cK}{\mathcal{K}}
\newcommand{\cC}{\mathcal{C}}
\newcommand{\cS}{\mathscr{S}}
\newcommand{\cP}{\mathcal{P}}
\newcommand{\cQ}{\mathcal{Q}}
\newcommand{\cD}{\partial}
\newcommand{\DD}{\mathsf{D}}
\newcommand{\cX}{\mathcal{X}}
\newcommand{\cx}{\mathbf{x}}
\newcommand{\cs}{\mathbf{s}}
\newcommand{\cc}{\mathbf{c}}
\newcommand{\cA}{\mathbf{A}}
\newcommand{\cb}{\mathbf{b}}
\newcommand{\cy}{\mathbf{y}}
\newcommand{\ccr}{\mathbf{r}}
\newcommand{\cu}{\mathbf{u}}
\newcommand{\cv}{\mathbf{v}}
\newcommand{\cQMat}{\mathbf{Q}}

\newcommand{\cz}{\mathbf{z}}
\newcommand{\cM}{\mathbf{M}}
\newcommand{\cg}{\mathbf{g}}
\newcommand{\ctheta}{\boldsymbol{\Theta}}
\newcommand{\ceta}{\boldsymbol{\eta}}
\newcommand{\cnu}{\boldsymbol{\nu}}
\newcommand{\cW}{\mathbf{W}}
\newcommand{\acC}{\mathbf{C}}

\newboolean{showcomments}
\setboolean{showcomments}{true}
\newcommand{\ja}[1]{\ifthenelse{\boolean{showcomments}}
{ \textcolor{red}{#1}}{}}
\newcommand{\fz}[1]{\ifthenelse{\boolean{showcomments}}
{ \textcolor{blue}{#1}}{}}

\newtheorem{remark}{Remark}

\newtheorem{proposition}{Proposition}
\newtheorem{definition}{Definition}
\newtheorem{corollary}{Corollary}
\newtheorem{lemma}{Lemma}
\newtheorem{theorem}{Theorem}

\title{The Optimal Power Flow Operator: Theory and
Computation \thanks{Submitted to the editors DATE.
\funding{This work was funded by NSF grants CCF 1637598, CPS 1739355, and ECCS 1619352, PNNL grant 424858, ARPA-E NODES through grant DE-AR0000699, and through the ARPA-E GRID DATA program.  }}}

\begin{document}
\title{The Optimal Power Flow Operator: Theory and Computation}


\author{Fengyu~Zhou,~\IEEEmembership{Student Member,~IEEE,}
        James~Anderson,~\IEEEmembership{Member,~IEEE,}
        and~Steven~H.~Low,~\IEEEmembership{Fellow,~IEEE}% <-this % stops a space
\thanks{This work was funded by NSF grants CCF 1637598, CPS 1739355, and ECCS 1619352, PNNL grant 424858, ARPA-E NODES through grant DE-AR0000699, and through the ARPA-E GRID DATA program.}
\thanks{Fengyu Zhou and Steven H. Low are with the Department of Electrical Engineering, California Institute of Technology, Pasadena,
CA 91125 USA (e-mail: \{f.zhou, slow\}@caltech.edu). }% <-this % stops a space
\thanks{James Anderson  is now with the Department of Electrical Engineering, Columbia University, New York, NY 10027 USA (e-mail: james.anderson@columbia.edu).}}


\markboth{}%
{Zhou \MakeLowercase{\textit{et al.}}: The Optimal Power Flow Operator}



\maketitle

\begin{abstract}
Optimal power flow problems (OPFs) are mathematical programs used to determine how to distribute power over networks subject to network operation constraints and the physics of power flows. In this work, we take the view of treating an OPF problem as an operator which maps user demand to generated power, and allow the network parameters (such as generator and power flow limits) to take values in some admissible set. The contributions of this paper are to formalize this operator theoretic approach,  define and characterize  restricted parameter sets under which the mapping has a singleton output, independent binding constraints, and is differentiable. In contrast to related results in the optimization literature, we do not rely on introducing auxiliary slack variables. Indeed,  our approach provides results that have a clear interpretation with respect to the power network under study. We further provide a closed-form expression for the Jacobian matrix of the OPF operator and describe how various derivatives can be computed using a recently proposed scheme based on homogenous self-dual embedding. Our framework of treating a mathematical program as an operator allows us to pose sensitivity and robustness questions from a completely different mathematical perspective and provide new insights into well studied problems. 
\end{abstract}

%Optimal power flow problems (OPFs) are mathematical programs used to distribute power over networks subject to network operation constraints and the physics of power flows. In this work we take the view of treating an OPF problem as an operator which maps user demand to generated power, and allow the network parameters to take values in some admissible set. The contributions of this paper are to formalize this operator theoretic approach,  define and characterize a restricted parameter sets under which the mapping has a singleton output, independent binding constraints, and is differentiable. We further provide a closed-form expression for the Jacobian matrix of the OPF operator and describe how various derivatives can be computed using a recently proposed scheme based on homogenous self-dual embedding. Our framework of treating a mathematical program as an operator allows us to pose sensitivity and robustness questions from a completely different mathematical perspective and provide new insights into well studied problems. 


% \fi
\IEEEpeerreviewmaketitle

\section{Introduction}

\IEEEPARstart{G}{iven} a power network, the optimal power flow (OPF) problem seeks to find an  operating point that minimizes an appropriate cost function subject to power flow constraints (e.g. Kirchhoff's laws) and pre-specified network tolerances (e.g. capacity constraints)\cite{Car62,HunG91,DomT68,FraSR12,Fra016}. The decision variables in an OPF problem are typically voltages and generation power. Cost function choices include minimizing power loss, generation cost, and user disutility.

In this work we consider the direct current (DC) model of the power flow equations~\cite{StoJA09dc, SunT10,Woo2012power}. The DC OPF problem is widely used in  industry and takes the form of a linear program (LP)~\cite{BerT97,Dan16}. 
%At a high level, for a network with $\nG$ generators and $\nL$ loads, the DC OPF problem takes the following form:
Consider a power network with $\nG$ generators and $\nL$ loads, the problem is formulated as
\begin{subequations}
\begin{eqnarray}
\underset{\sgen}{\text{minimize}}~~  && \f^{\T}\sgen
\label{eq:opf0.a}\\
\text{ subject to }& & \Aeq\sgen = \beq(\sload,\bpara),
\label{eq:opf0.b}\\
& &  \Ain\sgen \leq \bin
\end{eqnarray}
\label{eq:opf0}
\end{subequations}
where $\sgen \in \Real^{\nG}$ is the decision vector of power generations at each generator node in the network. The equality constraint function $\beq$ is linear in $\sload$ (the vector of power demands at each node) and $\bpara$. We view $\f$, $\bpara$ and $\bin$ as network parameters, whose values are allowed to change in a set $\Omega$. 
Matrices $\Aeq$ and $\Ain$ are determined by network topology and susceptances.
We treat the power demand $\sload$, which is allowed to take values in the set $\setsl$, as an ``input'' and the optimal generations $(\sgen)^{\star}$ as an ``output''. One of the main contributions of this work is to study the DC OPF problem
~\eqref{eq:opf0} as an operator, in particular we define
\begin{equation*}
\OPF :  \setsl \rightarrow 2^{\setS}
\end{equation*}
where $\setS \subset \Real^{\nG}$ and $2^{\setS}$ denotes the power set of $\setS$. 

%The first focus of this paper is to study if it is easy (and how easy) for the network parameters to endow $\OPF$ and the underlying OPF problem with all the following properties.
We begin by deriving conditions such that:
\begin{enumerate}
\item $\OPF$ maps to a singleton.
\item $\OPF$ is continuous everywhere and differentiable almost everywhere.
\item All points in the image space of $\OPF$ are optimal solutions to a DC-OPF problem with a fixed number of binding constraints (the number will be later derived as $\nG-1$).
%\item At the optimal point, OPF problem has independent and fixed number of binding constraints (the number will be later derived as $\nG-1$).
\end{enumerate}

Results such as those above have been shown to hold (with high probability) when the problem is simply viewed as a mathematical program (see for example;  \cite{pataki2001generic,dur2016genericity, luc1997differentiable,Fia1983introduction}), of which the DC-OPF problem is a special case. However these results, whilst insightful in optimization theory convey little actionable information about the process the optimization problem models.

We note  that Property 3 is equivalent to the concept of nondegeneracy \cite{wang1997nondegeneracy}. In contrast to  standard results,  we narrate from the perspective of binding constraints. Such a perspective is beneficial because it provides insight into the physical meaning of OPF problem solutions. In particular, each binding constraint implies that either a generator or power line is idle or saturated. Furthermore, we use a sequence of restrictions to obtain subsets of ``good'' parameters (i.e., those which preserve the properties listed), thus our method is both constructive and interpretable.

With  properties 1) -- 3) established, we provide a closed form expression for the Jacobian matrix of $\OPF$, that is, an approximate linear mapping from input (generator loads) to output (an optimal DC-OPF solution).  We also show how one could equivalently view the Jacobian  as a function of the set of binding constraints, i.e, a mapping which is independent of network parameters. This new perspective presents advantages especially if we are more interested in the global behavior of the Jacobian map (such as the worst-case amplification of a generator gain given a change in load), rather than just the generations given a specific load profile. 
Finally, we conclude by describing how various derivatives can be computed using  recently developed ideas~\cite{ArgBBBM19,AmoK17}.

From a power networks perspective, we are specifically interested in formalization since there has been a lot of research interest on the topic of characterizing the relationship between power load and generation over the last few years, see for example \cite{ZhoAL19,GenX16,ji2016probabilistic,misra2018learning, ng2018statistical},  and the references therein.

This work is a first step towards characterizing this complicated relationship. Specifically, establishing uniqueness of solution is a fundamental property of an operator as it provides the foundation for defining a derivative. Moreover, many numerical techniques require unique solutions to ensure convergence. 
Characterizing the set of  independent binding constraints paves the road for may further topics of interest. For instance, 
in \cite{RoaM19}, it shown that even under significant load variations, the number of binding line constraints in the DC-OPF problem is frequently a small proportion of the total number of constraints -- an observation which has significant implications when it comes to long-term planning and assessing network vulnerability~\cite{KimWBH16}. In \cite{GenX16}, the set of binding constraints determines an area of load profiles (termed System Pattern Regions  by the authors) within which the vector of locational marginal prices remain constant.

%\fz{Knowing OPF problem has independent binding constraints, the closed form expression for the Jacobian matrix can serve as the main tool to study the system sensitivity.}. %JA: This seems like a bit too uch information here

Sensitivity, more broadly is a useful term to quantify. In recent works we showed that  ``worst-case'' sensitivity bounds provide privacy guarantees when releasing power flow data~\cite{ZhoAL19} and for data disaggregation~\cite{AndZL18}. In the context of real-time optimization where sensitivity is often assumed to be known and bounded~\cite{TanDBL18}, this work can be used to provide exactly these bounds (or rule them out). 

\subsection{Related Work}
Our approach to sensitivity analysis differs from the standard perturbation approach which assumes the constraints are shifted from their nominal right-hand sides, and then looks at the Lagrange multipliers; see for example~\cite[Ch. 5.6]{Boy04convex}. In such a setting one considers the optimization problem 
\begin{equation}\label{eq:opt_pert}
\{\text{minimize}_\x~f_0(\x)~:~ f_i(\x) \le \mathbf{u}_i,~ \mathbf{h}_j(\mathbf{x}) = \mathbf{v}_j , \forall i,j\}.
\end{equation}
%\begin{align}\label{eq:opt_pert}
%\underset{\x}{\text{minimize}} \quad & \mathbf{f}_0(\mathbf{x}) \nonumber \\ 
%\text{subject to} \quad & \mathbf{f}_i(\mathbf{x}) \le \mathbf{u}_i, \quad i=1,\hdots,m, \\
%\quad & \mathbf{h}_j(\mathbf{x}) = \mathbf{v}_j, \quad j=1,\hdots,p. \nonumber
%\end{align}
The \emph{nominal} form of problem~\eqref{eq:opt_pert} has   $\mathbf{u}_i=0$ and $\mathbf{v}_j=0$  for all $i$ and $j$. Denote the nominal optimal value  by $\mathbf{x}^{\star}$. The \emph{perturbed} problem is obtained by adjusting the right-hand side of the equalities and inequalities, thus tightening or relaxing the constraints. The perturbed optimal value is denoted by $p^{\star}(\mathbf{u},\mathbf{v})$. Let $\lbold^{\star}$ and $\nubold^{\star}$ denote the vectors of optimal Lagrange multipliers for the nominal problem. Let us assume strong duality holds. Then we have the following well known lower-bound for the perturbed problem
$p^{\star}(\mathbf{u},\mathbf{v}) \ge p^{\star}(\mathbf{0},\mathbf{0}) - \lbold^{\star \T}\mathbf{u} - \nubold^{\star \T}\mathbf{u}.$ 
Additionally, the magnitude (and sign in the case of $\nubold^{\star}$) provide information regarding the sensitivity of $p^{\star}(\mathbf{u},\mathbf{v})$ with respect to constraints being tightened and relaxed. 
Concretely, the sensitivities are given by the relations $\lbold_i^{\star} =  -\partial_{\mathbf{u}_i}  p^{\star}(\mathbf{0},\mathbf{0})$ and $\nubold_j^{\star} =  -\partial_{\mathbf{v}_j}  p^{\star}(\mathbf{0},\mathbf{0}) $. 
This approach differs from our problem in the following ways. First, we focus on the perturbation of the optimal solution, rather than the optimal value. Second, we focus more on the set of  binding constraints and whether they are independent, as opposed to appealing to duality.

%Under differentiability assumptions on $p^{\star}(\mathbf{u},\mathbf{v})$, the optimal Lagrange multipliers characterize the local sensitivities of the optimal value with respect to constraint perturbations. 

A related body of work to compare our results to, is that of robust optimization~\cite{BenGN09, BerBC11, MulVZ95} and stochastic optimization~\cite{HeyS04, KalW94,MilW65}.  In both cases, the goal is to mitigate the effects of uncertainty. In contrast, our work seeks to determine how the optimal decision changes with respect to data perturbations. Our results can thus be considered complementary to the robust and stochastic optimization frameworks.

 %Consider the linear program $\{\text{minimize}_\x~ \mathbf {c}^\T\x~:~ \mathbf{Ax} \le \mathbf{b}, (\mathbf{A,b,c})\in \mathcal{U}\}$. In the robust optimization setting, the set $\mathcal{U}$ encodes deterministic uncertainty in the problem data. The objective is thus to minimize the decision vector $x$ over all allowable instantiations of the uncertainty. In the stochastic setting, the same philosophy is applied. The set $\mathcal U$ encodes distributions from which the problem data is drawn and the cost function is then suitably defined. 
 
Within the optimization community there has been a lot of interest in sensitivity problems. Although not stated formally, some work has studied LPs as operators.  We can trace these sensitivity-type results back to Hoffman~\cite{Hof03} who showed that for vector $\x$ which satisfies $\mathbf{Ax\le b}$ and  $\z$ such that $\mathbf{Az \nleq b}$, then $\mathrm{dist}(\x,\z)$ is bounded from above by $\alpha\|(\mathbf{Az-b})_+\|$, where the constant $\alpha$ depends on $\mathbf A$. The value of $\alpha$ depends on the choices of norms and is typically difficult to compute. Work by Robinson~\cite{Rob75} addressed the problem of determining if a set of inequalities (defined on a Banach space) remains solvable when the right hand side vector is perturbed, achieving  Hoffman's result as a special case. In~\cite{ManS87}, results based on perturbations to the right-hand-side of the inequalities in a linear program are shown to satisfy $\|\hat{\x} -\bar{\x}\|\le \beta \|\hat{\mathbf b}-\bar{\mathbf b} \|$, where  $\mathbf{A}\hat{\x}\le \hat{\mathbf b},  \mathbf{A}\bar{\x}\le \bar{\mathbf b}$, and the constant $\beta$ depends on $\mathbf{A}$. In some cases $\beta$ (the Lipschitz constant) can be found by solving a linear program. It is also shown that perturbations to the objective function destroy this Lipschitz continuity. A contribution of our work is to provide dense sets for which perturbations of the cost vector (and right-hand-side vectors) maintain differentiability.

With respect to uniqueness of solution, in \cite{pataki2001generic,dur2016genericity}, it has been proved that among all the linear program problem instances, almost all of those instances have unique (thus basic) and non-degenerate primal and dual solutions, and strict complementarity holds almost everywhere.
In the context of OPF problems, our results are more specific; First, traditional result shows nondegenerate instances are almost everywhere among all the problem instances, 
while our result says the good instances within a subset are also almost everywhere.
Taking \eqref{eq:opf0} as an example, traditional results show that for almost all the $(\f,\Aeq,\Ain,\beq,\bin)$, \eqref{eq:opf0} maintains the three desirable properties above.
In our work, we show that given any fixed $(\Aeq,\Ain)$ corresponding to power system structures, then for almost all the $(\f,\bpara,\bin)$, all the good properties listed above hold for almost all instances of $\sload$.
%{\bf We also want to point out that our results can be alternatively proved by extending or modifying some intermediate results in the literatures listed above,
%but we think deriving from basic linear algebra as we did in this paper is still the most concise and straightforward way.}
Second, traditional results are usually formulated in the canonical form in order to capture the general features of LPs. 
Our result, on the other hand, 
does not rely on introducing auxiliary slack variables and thus
reveals how those properties link to the physical behavior of power systems. 
%For example, having $\nG-1$ independent binding constraints in OPF directly means the total number of idle or saturated generators and power lines will always amount to exactly $1$ less than the total number of generators. This result turns to be useful in many applications \cite{ZhoAL19}.
%The differentiability for linear and nonlinear programs is well studied in many literatures \cite{luc1997differentiable,Fia1983introduction}.
%We leveraged those existing result to show the differentiability for OPF operator.
%We found that many core properties of LP has been well studied in theoretical optimization community, and they turn out to be very important if we want to formalize the OPF operator and further study the relationship between load and the optimal generation.
%To the best of our knowledge, few research works have rigorously studied and interpreted those properties in the context of power systems.
%We hope this work can abridge this gap.

There has been some work which  specifically defines and studies the DC-OPF sensitivity. In \cite{Gri1990optimal,Yu2001sensitivity}, the OPF problem is formulated as a parameterized optimization problem~\cite{GudVJ90} where the loads, the upper and lower bounds for generations and branch power flows are all parameterized by a single parameter. Under the assumptions that the binding constraints are known and the optimal solution and Lagrange multipliers are available, sensitivity can be computed. The method is restrictive because there is only a single degree of freedom in  parameter variation. Moreover, determining differentiability is a very involved process. On the contrary, we generalize the concept of sensitivity to the Jacobian matrix, which allows the parameters to change in various directions.
Instead of checking differentiability for each problem,
we explicitly characterize the sets of parameters that guarantee differentiability, and based on the fact that those sets are all dense within the spaces of interest, we conclude that differentiability can always be assumed up to parameter perturbation. Furthermore, we provide numerical methods to compute the derivatives.

\subsection{Article Outline}
In Section \ref{sec:background}, we formalize the DC OPF problem, and characterize the parameter set of interest under which the OPF problem is feasible and has optimal solutions.
For the parameters in the set of interest, we define the associated operator that maps the load to the set of optimal generations.
In Section \ref{sec:property}, we restrict the set of interest to those parameters that endow the OPF operator with desirable properties. We show that the restricted set is dense within the set of interest, so the restriction does not lose generality up to perturbation.
In Section \ref{sec:derivative}, we prove the operator is differentiable and derive the closed form expression of the Jacobian matrix in terms of the independent binding constraints.
Moreover, we prove there exist a surjection between the restricted set and the set of independent binding constraints such that the derivative of the operator and the Jacobian matrix in terms of binding constraints take the same value under such surjection.
%While the former characterizes how sensitive the OPF operator is, the latter only depends on the graph topology and has a closed form expression.
In Section \ref{sec:comp}, we  demonstrate that an algorithm introduced in recent work \cite{ArgBBBM19,AmoK17} can help numerically evaluate the operator differentiation and extend the results to the alternating current (AC) OPF case.
Finally, Section \ref{sec:example} provides an illustrative example.%The DC-OPF problem, as well as the corresponding operator are formalized in \cref{sec:background}.
\section{Background}\label{sec:background}
%In this section we define the power network model and the optimal power flow problem. 
%We introduce some assumptions on allowable parameter sets and show that these assumptions are mild.

\subsection*{Notation}
Vectors and matrices are typically written in bold while scalars are not. Given two vectors $\mathbf{a},\mathbf{b} \in \Real^n$, $\mathbf a\ge \mathbf b$ denotes the element-wise partial order $\mathbf{a}_i \ge \mathbf{b}_i$ for $i=1,\hdots,n$. For a scalar $k$, we define the projection operator $[k]^-:= \min \{0,k\}$. We define $\|\mathbf x\|_0$ as the number of non-zero elements of the vector $\mathbf x$. 
Identity and zero matrices are denoted by $\eye^n$ and $\zero^{n\times m}$ while vectors of all ones are denoted by $\one_n$ where superscripts and subscripts indicate their dimensions.
To streamline notation, we omit the dimensions when the context makes it clear. 
The notation $\Real_+$ denotes the non-negative real set $[0,+\infty)$. %, and similarly $\Real_-$ denotes the nonpositive real set $(-\infty,0]$.
For $\mathbf X \in \Real^{n\times m}$, the restriction $\mathbf{X}_{\{1,3,5\}}$ denotes the $3\times m$  matrix composed of stacking rows $1,3$, and $5$ on top of each other. We will frequently use a set to describe the rows we wish to form the restriction from, in this case we assume the elements of the set are arranged in increasing order.
We will use $\base_{\scriptscriptstyle m}$ to denote the standard base for the $m$\textsuperscript{th} coordinate,  its dimension will be clear from the context. 
Let $(\cdot)^{\dagger}$ be the Moore-Penrose inverse. 
Denote $[m]:=\{1,2,\dots,m\}$ and $[n,m]:=\{n,n+1,\dots,m\}$.
Finally, for a convex set $\cX\subseteq\Real^{n}$ and vector $\cx\in\Real^{n}$, we let $\cP_{\cX}\cx$ be the projection of $\cx$ onto the set $\cX$. By isometry, the domain of the projection operator is extended to matrices when needed. 

\subsection{System Model}
Consider a power network modeled by an undirected connected graph $\Graph(\Vertex, \Edge)$, where $\Vertex:=\VertexG\cup\VertexL$ denotes the set of buses which can be further classified into subsets of generators $\VertexG$ and loads $\VertexL$, and $\Edge\subseteq \Vertex\times\Vertex$ is the set of all branches linking those buses. We will later use the terms (graph, vertex, edge) and (power network, bus, branch) interchangeably. Suppose $\VertexG\cap\VertexL=\emptyset$ and there are $|\VertexG|=:\nG$ generator and $|\VertexL|=:\nL$ loads, respectively.
For simplicity, let $\VertexG=[\nG]$, $\VertexL=[\nG+1,\nG+\nL]$.  Let $\nGL=\nG+\nL$. Without loss of generality, $\Graph$ is a connected graph with $|\Edge|=:\nE$ edges labelled as $1,2,\dots,\nE$.
Let $\CM\in\Real^{\nGL\times\nE}$ be the  incidence matrix, where the orientation of the edge is arbitrarily chosen.   We will use $e$, $(u,v)$ or $(v,u)$ interchangeably to denote the same edge.
Let $\B={\rm diag}(b_1,b_2,\dots,b_E)$, where $b_e>0$ is the susceptance of branch $e$. As we adopt a
DC power flow model, all branches are assumed lossless. Further, we denote the generation and load as $\sgen\in\Real^{\nG}$, $\sload\in\Real^{\nL}$, respectively.
Thus $\sgen_i$ refers to the generation on bus $i$ while $\sload_i$ refers to the load on bus $\nG+i$. We will refer to bus $\nG+i$ simply as load $i$ for simplicity. 
The power flow on branch $e\in\Edge$ is denoted as $\pf_e$, and $\pf:=[\pf_1,\dots,\pf_{\nE}]^{\T}\in\Real^\nE$ is the vector of all branch power flows.
To simplify analysis, we assume that there are no buses in the network that are both loads and generators. This is stated formally below:
\begin{assumption}\label{A1}
$\VertexG\cap\VertexL=\emptyset$.%  and $\deg(i)=1$ for all $i \in \VertexG$.
\end{assumption}

The above assumption is not restrictive.
We can always split a bus with both a generator and a load into a bus with
only the generator adjacent to another bus with only the load, and
connect all the neighbors of the original bus to that load bus. % as shown in \cref{fig:A1}.

%\begin{figure}
%\centering
%\includegraphics[width=0.45\columnwidth]{./Figures/Fig_for_A1.pdf}
%\caption{Suppose bus $i$ is a generator with load, and is adjacent to 3 other buses, then we can create a new network where $i$ is replaced by $i'$, a pure generator, and $i''$, a pure load. The original neighbors of $i$ are adjacent to $i''$ in the new network. The susceptance for branch $(i',i'')$ can be assigned with any positive value. Two networks above are equivalent if we view $i'$ and $i''$ as a whole.}
%\label{fig:A1}
%\end{figure}

\subsection{DC Optimal Power Flow}\label{sec:OPF}
We focus on the DC-OPF problem with a linear cost function \cite{StoJA09dc, SunT10,Woo2012power}. That is to say, the voltage magnitudes are assumed to be fixed and known and the lines are considered to be lossless. Without loss of generality, we assume all the voltage magnitudes to be $1$. The decision variables are the voltage angles denoted by vector $\ang\in\Real^{\nGL}$ and power generations $\sgen$, given loads $\sload$. The DC-OPF problem takes the form:
\begin{subequations}
\begin{eqnarray}
\underset{\sgen, \ang}{\text{minimize}} ~~ && \f^{\T}\sgen
\label{eq:opf1.a}\\
\text{ subject to }& & \ang_1 = 0
\label{eq:opf1.b}\\
& &  \CM\B\CM^{\T} \ang = 
\left[ \begin{array}{c}
\sgen \\
-\sload
\end{array} \right] 
\label{eq:opf1.c}\\
& & \genllim \leq\sgen\leq \genulim
\label{eq:opf1.d}\\
& &  \underline{\pflow}\leq\B\CM^{\T}\ang\leq\overline{\pflow}.
\label{eq:opf1.e}
\end{eqnarray}
\label{eq:opf1}
\end{subequations}
Here, $\f\in\Real_+^{\nG}$ is the unit cost for each generator, and bus $1$ is selected as the slack bus with fixed voltage angle $0$. In \eqref{eq:opf1.c}, we define the injections for generators to be positive, while  injections for loads are defined as the negation of $\sload$.
The upper and lower limits on the generations are set as $\genulim$ and $\genllim$, respectively, and $\overline{\pflow}$ and $\underline{\pflow}$ are the limits on branch power flows. 
We assume that~\eqref{eq:opf1} is well posed, i.e. $\genulim>\genllim \ge 0$, $\overline{\pflow}>\underline{\pflow}$. Note that the LP~\eqref{eq:opf1} is a particular realization of~\eqref{eq:opf0}.\footnote{Though two problems have different decision variables, one can always replace $\ang$ in \eqref{eq:opf1}  by $(\CM\B\CM^\T)^{\dagger}[(\sgen)^\T,(-\sload)^\T]^\T$ to absorb the additional decision variable $\ang$.}
%Our results can be extended to the general case, however we present this paper assuming the specific problem form of ~\eqref{eq:opf1}. We stress that no knowledge of power engineering is needed to derive or understand the results in this paper.

Let $\taueq\in\Real^{\nGL+1}$ be the vector of Lagrangian multipliers 
associated with equality constraints \eqref{eq:opf1.b}, \eqref{eq:opf1.c},
and $(\lambdap,\lambdam)$ and $(\mup,\mum)$ be the Lagrangian multipliers 
associated with inequalities \eqref{eq:opf1.d} and \eqref{eq:opf1.e} respectively. 
As \eqref{eq:opf1} is a linear program, the following KKT condition
holds at an optimal point when \eqref{eq:opf1} is feasible:
\begin{subequations}
\begin{eqnarray}
&& \eqref{eq:opf1.b}-\eqref{eq:opf1.e}
\label{eq:KKT.a}\\
&& \bf{0}=M^{\T}\taueq+\CM\B(\mup-\mum)
\label{eq:KKT.b}\\
&&  -\f=-[\taueq_1,\taueq_2,\cdots,\taueq_{\nG}]^{\T}+\lambdap-\lambdam
\label{eq:KKT.c}\\
&& \mup,\mum,\lambdap,\lambdam\geq 0
\label{eq:KKT.d}\\
&& \mup^{\T}(\B\CM^{\T}\ang-\overline{\pflow})=\mum^{\T}(\underline{\pflow}-\B\CM^{\T}\ang)=0
\label{eq:KKT.e}\\
&& \lambdap^{\T}(\sgen-\genulim)=\lambdam^{\T}(\genllim-\sgen)=0,
\label{eq:KKT.f}
\end{eqnarray}
\label{eq:KKT}
\end{subequations}
where
\begin{align*}
\M:=\left[
\begin{array}{c}
\CM\B\CM^{\T}\\
 \base_1^{\T}
\end{array} 
\right]
\end{align*}
is an $(\nGL+1)$-by-$\nGL$ matrix with rank $\nGL$, and $\base_1$ denotes the standard first basis vector. 
Condition \eqref{eq:KKT.a} corresponds to primal feasibility, condition \eqref{eq:KKT.d} corresponds to dual feasibility, conditions \eqref{eq:KKT.e}, \eqref{eq:KKT.f} correspond to complementary slackness,
and conditions \eqref{eq:KKT.b}, \eqref{eq:KKT.c} correspond to stationarity.

\subsection{OPF as an Operator: $\OPF$}\label{sec:assump}
We will now describe how to formulate the DC-OPF~\eqref{eq:opf1} as a mapping from load to (optimal) generation space. We assume throughout the paper that the topology of the network remains constant, as do the line susceptances. These assumptions imply that the graph Laplacian given by $\CM\B\CM^{\T} $ does not change.
Let $\para:=[(\genulim)^{\T},(\genllim)^{\T},\overline{\pflow}^{\T},\underline{\pflow}^{\T}]^{\T}\in\Real^{2\nG+2\nE}$ be the vector of system limits. Define
\begin{align}\nonumber
\setpara:=\{\para~|~\genllim\geq 0, \eqref{eq:opf1.b}-\eqref{eq:opf1.e}~\text{are feasible for some}~\sload>0\}.
\end{align}
The set $\setpara$ defines the set of power flow and generation limits such that the DC-OPF is primal-dual feasible and makes physical sense i.e. upper-limits are greater than lower-limits. %Note that $\setpara$ does not depend on the cost vector $\f$.
%\begin{itemize}
%\item $\exists\sload>0$, such that \ja{the DC-OPF constraints} \eqref{eq:opf1.b}-\eqref{eq:opf1.e} are feasible;
%\item $\genllim\geq 0$.
%\end{itemize}

%\footnote{This condition will be satisfied if 
%\begin{eqnarray*}
%\setpara & := & \prod\Omega_{\genulim_i}\times\prod\Omega_{\genllim_i}\times\prod\Omega_{{\overline{\pflow}}_i}\times\prod\Omega_{{\underline{\pflow}}_i}
%\end{eqnarray*}
%where $\Omega_{\genulim_i}$, $\Omega_{\genllim_i}$, $\Omega_{{\overline{\pflow}}_i}$ and $\Omega_{{\underline{\pflow}}_i}$ are the set of limits $\genulim_i, \genllim_i, \overline{\pflow}_i,
%\underline{\pflow}_i$ respectively.  This assumes that these limits are
%independent of each other.}
For each $\para\in\setpara$, define\footnote{In practice, if a load has $0$ value, one could replace it by an arbitrarily small positive value so that the load profile is always strictly positive.}
\begin{align}\nonumber
\setsl(\para):=\{\sload~|~\sload>0, \eqref{eq:opf1.b}-\eqref{eq:opf1.e}~\text{are feasible}\}.
\end{align}
% $\setsl(\para)$ be the corresponding set of $\sload$ such that 
%\begin{itemize}
%\item constraints \eqref{eq:opf1.b}-\eqref{eq:opf1.e} are feasible.
%\item $\sload>0$
%\footnote{Since the feasible domain of \eqref{eq:opf1} is compact, the optimal solution always exists.}
% \item other potential constraints on $\sload$ implied by the realistic system.
%\end{itemize}
Then $\setsl(\para)$ is convex and nonempty.
%\footnote{This condition will be satisfied if $\setsl:=\prod\Omega_{\sload_i}$,
%where $\Omega_{\sload_i}$ are the set of $\sload_i$. This assumes that these loads are
%independent of each other.}
%In practice, the loads are always within their own ranges $\Omega_{\sload_i}$, which are also independent of each other, so one commonly used example for $\setsl$ is $\prod\Omega_{\sload_i}$, which satisfies the assumption.
When we fix $\para$ and there is no confusion, we simply write $\setsl$.
\begin{definition}\label{df:setall}
Define $\setall:=\{(\f,\para,\sload)~|~\f\in\Real_+^{\nG},\para\in\setpara,\sload\in\setsl(\para)\}$.
\end{definition}
When $\para\in\setpara$ and $\sload\in\setsl(\para)$ the DC-OPF problem \eqref{eq:opf1} is feasible.
As \eqref{eq:opf1.b} fixes the angle $\ang_1$ at the slack bus, and \eqref{eq:opf1.b} restricts the angle difference between any two adjacent buses we conclude that the feasible set of \eqref{eq:opf1} is compact, and thus, by Weirstrass' Theorem, the optimal solutions to \eqref{eq:opf1} always exist. 
We now define the operator $\OPF$, which is the central object of study in this paper.

\begin{definition}
Fix $\para\in\setpara$ and $\sload\in\setsl(\para)$,
let the set valued operator $\OPF:\setsl\rightarrow 2^{\Real^{\nG}}$ be the mapping such that $\OPF(\bf{x})$ is the set of optimal solutions to \eqref{eq:opf1} with parameter $\sload=\bf{x}$.
\end{definition}
In the following section we will establish various properties of the $\OPF$ operator and show that it is a valuable tool for gaining insight into the sensitivity, robustness, and structure of the DC OPF problem~\eqref{eq:opf1}.
\section{Operator Properties}\label{sec:property}
We assume that the network topology and line susceptance are fixed, that is  $\CM$ and $\B$ are constant.
The operator $\OPF$ is parameterized by  $\f$, $\para$ and $\sload$.
The set $\setall$ defined in Definition \ref{df:setall} prescribes all the parameters under which \eqref{eq:opf1} is feasible.
%In this section, we restrict the parameter set so as to endow the operator $\OPF$ with desirable properties.
\subsection{Uniqueness}
We are specifically interested in the case when the OPF operator defined above maps to a singleton.
To pave the way for further properties in the following subsections, we in fact consider the vector $\f$  under heavier constraints.
Let $\setf$ be the set of vectors $\f\geq0$ such that $\forall \para\in\setpara, \forall\sload\in\setsl(\para)$:
\begin{itemize}
\item DC-OPF problem \eqref{eq:opf1} has a unique solution.
\item The KKT-system \eqref{eq:KKT} satisfies
\begin{align}\label{eq:multiplier1}
\|\mup\|_0+\|\mum\|_0+\|\lambdap\|_0+\|\lambdam\|_0\geq \nG-1.
\end{align}
\end{itemize}
\begin{proposition}\label{fisdense}
$\setf$ is dense in $\Real_+^{\nG}$.
\end{proposition}
\begin{IEEEproof}
See Appendix \ref{app:proof_fisdense}. 
\end{IEEEproof}

Proposition \ref{fisdense} shows that for a fixed network, it is easy to find an objective vector $\f$ such that \eqref{eq:opf1} not only has a unique solution for feasible $\sload$,
but also gives sufficiently many non-zero dual variables. 
%Thus we can make the following assumption and define the OPF operator accordingly:
For the remainder of the paper, the following assumption is in play:
\begin{assumption}\label{A:vectorf}
The objective vector $\f$ is in $\setf$.% i.e, $\f$  always guarantees the uniqueness of the solution to \eqref{eq:opf1} for all $\sload\in\setsl$.
\end{assumption}
This assumption ensures that~\eqref{eq:opf1} has to have a unique solution. When Assumption \ref{A:vectorf} does not hold, Proposition \ref{fisdense} implies that we can always perturb $\f$ such that the assumption is valid.
\begin{remark}
Under Assumption \ref{A:vectorf}, the value of $\OPF$ is always a singleton, 
so we can overload $\OPF({\bf x})$ as the function mapping from ${\bf x}$ to the unique optimal solution of \eqref{eq:opf1} with parameter $\sload=\bf{x}$.\footnote{Except for Appendix \ref{app:proof_wc_reduction} where $\OPF$ is still viewed as a set valued function, $\OPF$ will be viewed as a vector valued function throughout the paper by default.}
Since the solution set to the parametric linear program is both upper and lower hemi-continuous \cite{Zha1990note}, $\OPF$ is continuous as well.
\end{remark}

\subsection{Independent Binding Constraints}
The analysis on the OPF operator can usually be simplified if the set of binding (active) constraints at the optimal point is independent. Here, binding constraints refer to the set of equality constraints \eqref{eq:opf1.b}, \eqref{eq:opf1.c},
 \emph{and}  those inequality constraints \eqref{eq:opf1.d}, \eqref{eq:opf1.e}
  for which either the  upper or lower-bounds are active. Grouping the coefficients of these constraints into a single matrix $\bf{Z}$ we refer to them as being independent if $\bf{Z}$ is full-rank. Finally, define
%Let $\setslr(\para)$ be the set
\begin{align*}
\setslr(\para,\f):=\{&\sload\in\setsl(\para)~|~\eqref{eq:opf1}~\text{has exactly}~\nG-1\text{ binding}\\
&\text{inequalities }
\text{at the optimal point, given }\sload\}.
\end{align*}
When $\f$ is fixed, we shorten $\setslr(\para,\f)$ as $\setslr(\para)$. Further, if $\para$ is also fixed, then we will simply use $\setslr$.
%Then we have the following proposition.
\begin{theorem}\label{limitsdense}
For a fixed $\f\in\setf$, there exists a dense set $\setparar(\f)\subseteq\setpara$ such that $\forall \para\in\setparar(\f)$, 
the following statements are true:
\begin{itemize}
\item $\clos(\inte(\setsl(\para)))=\clos(\setsl(\para))$.
\item $\setslr(\para,\f)$  is dense in $\setsl(\para)$.
\end{itemize}
\end{theorem}
\begin{IEEEproof}
See Appendix \ref{app:proof_limitsdense}.
\end{IEEEproof}

%Finally,  \cref{A:limits} restricts the value of limits in \eqref{eq:opf1} in order to further simplify our analysis.
\begin{assumption}\label{A:limits}
The parameter $\para$ for the limits of generations and branch power flows is assumed to be in $\setparar(\f)$, as proposed in Theorem \ref{limitsdense}. 
\end{assumption}

Assumption \ref{A:limits} allows one to work with sets $\setparar(\f)$ that are well behaved (where ``well behaved'' is interpreted as $\setsl$ and $\inte(\setsl)$ having the same closure
%\ja{[solid sets? I can't find a definition for this in any topology text book. O know what a solid cone is, what is a solid set? A set with a non-empty relative interior?]}
and there being exactly $N_G-1$ binding constraints at the optimal point in the associated DC-OPF problem for almost every $\sload$). This assumption is important as in Section \ref{sec:derivative} it will be needed to show that the derivative of $\OPF$ exists almost everywhere. If Assumption \ref{A:limits} does not hold, Theorem \ref{limitsdense} implies that we can always perturb $\para$ such that the assumption holds. 
In the context of DC OPF problem, it also means for almost all the problem instances, there are exactly $\nG-1$ binding constraints at the optimal point.
The proof of Theorem \ref{limitsdense} can directly extend to the following two corollaries:
%Proofs are omitted due to the limited space.
\begin{corollary}\label{Co:cover_by_planes}
 The set $\setsl\setminus\setslr$ can be covered by the union of finitely many affine hyperplanes.
\end{corollary}
\begin{corollary}\label{Co:independent}
For any $\sload\in\setslr$, the $\nG-1$ tight inequalities in \eqref{eq:opf1}, along with $\nGL+1$ equality constraints, are independent.
\end{corollary}
\begin{definition}
Define $\setallr:=\{(\f,\para,\sload)~|~\f\in\setf,\para\in\setparar(\f),\sload\in\setslr(\para,\f)\}$.
\end{definition}

\begin{figure}
\centering
\includegraphics[width=1.0\columnwidth]{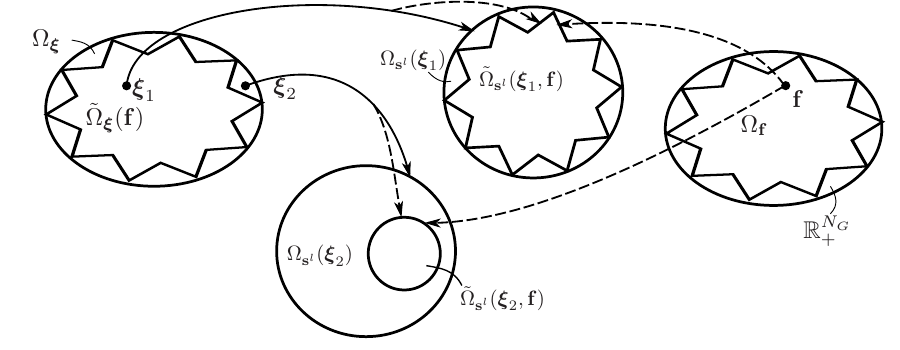}
\caption{Relationship among definitions in Section \ref{sec:assump}. Solid arrows show the mapping from $\para$ to $\setsl(\para)$, and dashed arrows show the mapping from $(\para,\f)$ to $\setslr(\para,\f)$. 
A star inscribed within an oval indicates the former set is dense within the latter.}
\label{Fig:Defs2}
\end{figure}

In summary, the two sets $\Omega$ and $\setallr$ characterize sets of objective functions, network parameters, and ``inputs'' that endow $\OPF$ with desirable properties. In particular  $\setall$ guarantees \eqref{eq:opf1} is feasible and $\OPF$ is thereby well-defined. 
The parameters in $\setallr$ additionally guarantee that \eqref{eq:opf1} has independent binding constraints and $\OPF$ is singleton-valued, and as will be shown in the next section, $\OPF$ is differentiable when $(\f,\para,\sload)\in\setallr$. The relationship among the sets $\setf$, $\setsl$, $\setslr$, $\setpara$, $\setparar$ defined above is illustrated in  Figure~\ref{Fig:Defs2}. Recall that 
informally, the set $\setpara$ contains all the $\para$ that make the OPF problem feasible, and $\setf$ contains $\f$ that guarantee the unique optimal solution for feasible OPF problems and sufficiently many non-zero Lagrange multipliers.
Proposition \ref{fisdense} shows $\setf$ is dense in $\Real_+^{\nG}$.
Each $\para\in\setpara$ maps to a set $\setsl(\para)$, while each $(\para,\f)$ maps to set $\setslr(\para,\f)$, which is a subset of $\setsl(\para)$. 
For fixed $\f$, by collecting all the $\para$ such that $\setsl(\para)$ has ``good'' topological property and $\setslr(\para,\f)$ is dense in $\setsl(\para)$, we obtain a set $\setparar(\f)$ depending on $\f$, 
and Proposition \ref{limitsdense} implies $\setparar(\f)$ is always dense in $\setpara$.

Since the sets that imply ``good'' properties ($\setf$, $\setslr$, $\setparar$) are all dense with respect to the corresponding whole sets of interest ($\Real_+^{\nG}$, $\setsl$, $\setpara$),
one can always perturb the parameters to endow $\OPF$  with these desirable properties.

\section{On the OPF Derivative}\label{sec:derivative}
In this section we show that  $\OPF$ is differentiable almost everywhere. 
We also provide an equivalent perspective from which to view the derivative (Jacobian matrix) of $\OPF$ in terms of binding constraints, and derive its closed form expression.\footnote{A word on notation is in order here. We denote the derivative of $f(x)$ with respect to $x$ by $\partial_x f$, however in some cases when there are complex dependencies on $x$ we will use $\frac{\partial f}{\partial x}$. In Section~\ref{sec:comp} when we deal with derivatives of conic programs we use the notationally lighter differential operator $\DD$.}
%\ja{\textbf{We need some text here...}}
\subsection{Existence}
%\ja{\textbf{[Notation in Sec 4.1 needs to be changed so that it's consistent with the rest of the paper: $x\rightarrow \x$]}}
Before deriving the expressions for the $\OPF$ derivative, it is necessary to guarantee that the operator is in fact differentiable.
The following lemma proposed in \cite{Fia1976sensitivity} and \cite{Fia1983introduction} gives the sufficient condition of differentiability. % for parameterized optimization problems.
We rephrase the lemma as follows. 
\begin{lemma}[\cite{Fia1976sensitivity,Fia1983introduction}]\label{Le:existence}
Consider a generic optimization problem parametrized by $\ctheta$:
\begin{subequations}
\begin{eqnarray}
\underset{\x\in\Real^n}{\rm{minimize}}  && ~f(\x;\ctheta)
\label{eq:opt.a}\\
\text{\rm  subject to }& & ~g_i(\x;\ctheta)\leq 0, i=1,2,\dots,m
\label{eq:opt.b}\\
& & ~h_j(\x;\ctheta)= 0, j=1,2,\dots,l.
\label{eq:opt.c}
\end{eqnarray}
\label{eq:opt}
\end{subequations}
If $(\x^*,\ceta^*,\cnu^*)$ is the primal-dual optimal solution for some $\ctheta_0$ and satisfies:
\begin{enumerate}[1)]
\item \label{it:c1}$\x^*$ is a locally unique primal solution.
\item \label{it:c2}$f,g_i,h_j$ are twice continuously differentiable in $\x$ and differentiable in $\ctheta$.
\item \label{it:c3}The gradients $\nabla g_i(\x^*)$ for binding inequality constraints and $\nabla h_j(\x^*)$ for equality constraints are independent.
\item \label{it:c4}Strict complementary slackness holds, i.e., $g_i(\x^*)=0\Rightarrow \ceta_i>0$.
\end{enumerate}
Then the local derivative $\partial_{\ctheta}\x^*$ exists at $\ctheta_0$, and the set of binding constraints is unchanged in a small neighborhood of $\ctheta_0$.
\end{lemma}

Using the set definitions from the previous section and the above lemma, we obtain the following result:
\begin{theorem}\label{thm:derivative_exist}
Under Assumptions \ref{A:vectorf} and \ref{A:limits}, for $\sload\in\setslr$, the derivative $\partial_{\sload}\OPF(\sload)$ always exists, 
and the set of binding constraints stay unchanged in some neighborhood of $\sload$.
\end{theorem}
\begin{IEEEproof}
By checking the conditions \ref{it:c1}-\ref{it:c4} in Lemma \ref{Le:existence}, the proof is established.
Alternatively, the theorem can be proved by extending Proposition 3.2 in \cite{luc1997differentiable}.
\end{IEEEproof}

Having established existence of the derivative of $\OPF$ we are now ready to study the associated Jacobian matrix.
\subsection{Jacobian Matrix}\label{subsec:Jacobian}
The Jacobian is an important tool in sensitivity analysis as it provides the best linear approximation of an operator from input to output space. The results of the previous section ensure that the partial derivatives exist almost everywhere. Let
\begin{align}\label{eq:Jacobian}
\JM(\sload;\f,\para):=\partial_{\sload}\OPF(\sload)
\end{align}
for $(\f,\para,\sload)\in\setallr$ denote the Jacobian of $\OPF$ at $\sload$.
To reduce the notational burden, we will simply use $\JM$ or $\JM(\sload)$ for short when the value of $(\f,\para,\sload)$ or $(\f,\para)$ is clear from  context.
Suppose at point $\sload$, the set of generators corresponding to binding inequalities is $\setSgen\subseteq\VertexG$, while the set of branches corresponding to binding inequalities is $\setSbra\subseteq\Edge$. From Theorem \ref{limitsdense} and Assumption \ref{A:limits} the following corollary is immediate:
\begin{corollary} $|\setSgen|+|\setSbra|=\nG-1$.
\end{corollary}
As Lemma \ref{Le:existence} implies that generators $\setSgen$ and branches $\setSbra$ still correspond to binding constraints near $\sload$, there is a local relationship between $(\sgen)^{\star}=\OPF(\sload)$ and $\sload$:
\begin{align}\label{eq:relationship}
\RelationMat
\left[
\begin{array}{c}
\sgen\\
\hline
\ang
\end{array}
\right]=
\left[
\begin{array}{c}
\zero^{\nG\times1}\\
\hline
-\sload\\
\gamma^{\T}\para\\
0
\end{array}
\right],
\RelationMat:=\left[
\begin{array}{c|c}
-\eye^{\nG} & \eye_{\VertexG}^{\nGL}\CM\B\CM^{\T}\\
\hline
\zero^{\nL\times\nG} & \eye_{\VertexL}^{\nGL}\CM\B\CM^{\T}\\
\eye_{\setSgen}^{\nG} & \zero^{|\setSgen|\times\nGL}\\
\zero^{|\setSbra|\times\nG} & \eye_{\setSbra}^{\nE}\B\CM^{\T}\\
\zero^{1\times\nG} & \base_1^{\T}
\end{array}
\right].
\end{align}
When there is no danger of confusion, we  use $\sgen$ to denote $(\sgen)^{\star}$.
\begin{figure}
\centering
\includegraphics[width=\columnwidth]{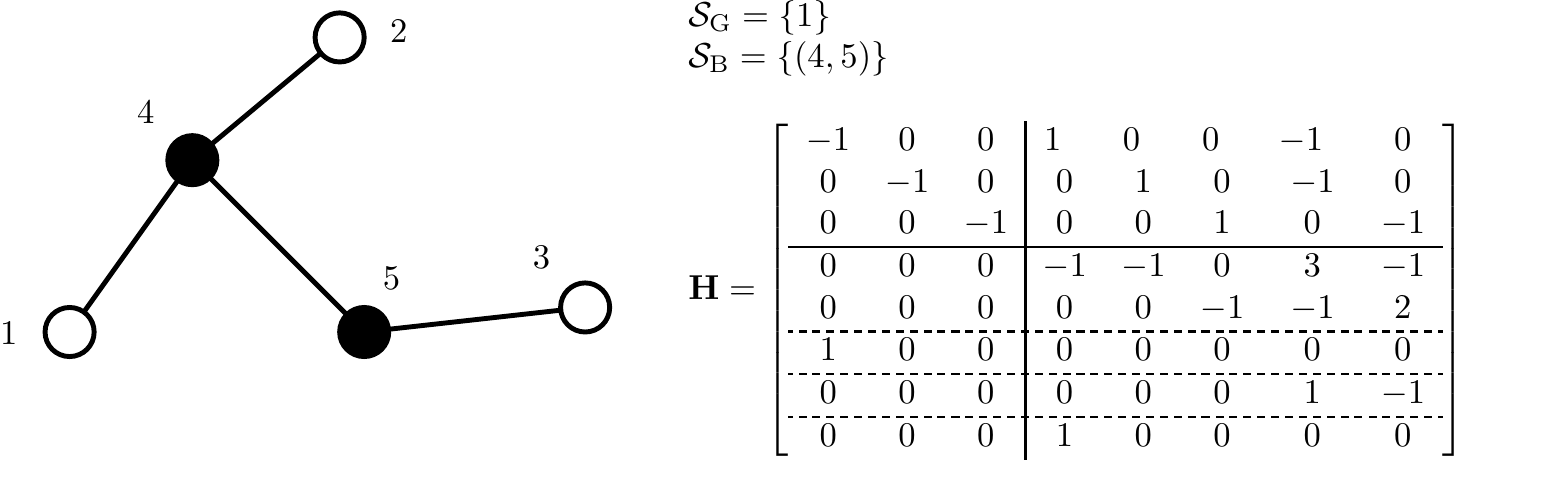}
\caption{In the above 5-bus network, white and black nodes denote generators and loads, respectively. Assume generator $1$ is a binding generator and branch $(4,5)$ is a binding branch. The 
$\RelationMat$ matrix for this example is given in the figure.}
\label{Fig:H}
\end{figure}
An example of $\RelationMat$ is given in Fig. \ref{Fig:H}.
On the right hand side, $\gamma\in\Real^{(2\nG+2\nE)\times(\nG-1)}$
where each column of $\gamma$ is a basis vector such that
$\gamma^{\T}\para$ gives a vector of capacity limits that binding generations and branch power flows hit. By Corollary \ref{Co:independent}, the first $\nGL+\nG-1$ rows of $\RelationMat$ are independent, and clearly the last row $[\zero,\base_1^{\T}]$ does not depend on the first $\nGL+\nG-1$ rows.
Hence $\RelationMat$ is invertible, and using the block matrix inversion formula, we have
\begin{align}
\left[
\begin{array}{c}
\sgen\\
\hline
\ang
\end{array}
\right]&=\RelationMat^{-1}
\left[
\begin{array}{c}
\zero^{\nG\times1}\\
\hline
-\sload\\
\gamma^{\T}\para\\
0
\end{array}
\right]
=\left[
\begin{array}{c|c}
* & \subRelationMat\\
\hline
* & *
\end{array}
\right]
\left[
\begin{array}{c}
\zero^{\nG\times1}\\
\hline
-\sload\\
\gamma^{\T}\para\\
0
\end{array}
\right]
\label{eq:relationinverse}
\end{align}
with  $\subRelationMat= \eye_{\tiny\substack{\VertexG}}^{\nGL}\CM\B\CM^{\T}\big(\bigMat(\setSgen,\setSbra)^{\T}\big)^{-1}$ and 
\begin{align}\label{eq:subMat}
%\nonumber
%\subRelationMat=&\eye_{\tiny\substack{\VertexG}}^{\nGL}\CM\B\CM^{\T}\big(\bigMat(\setSgen,\setSbra)^{\T}\big)^{-1},\\
\bigMat(\setSgen,\setSbra)^{\T}:=
&\left[\begin{array}{c}
 \eye_{\tiny\substack{\VertexL}}^{\nGL}\CM\B\CM^{\T}\\
\eye_{\setSgen}^{\nGL}\CM\B\CM^{\T}\\
\eye_{\setSbra}^{\nE} \B\CM^{\T}\\
\base_1^{\T}
\end{array}
\right].
\end{align}
Recall \eqref{eq:Jacobian} that $\sgen = \OPF(\sload)$ in \eqref{eq:relationinverse},
so the Jacobian matrix $\JM$ is
\begin{align}\label{eq:J}
\JM=-\subRelationMat (\eye_{[\nL]}^\nGL)^{\T}.
\end{align}
%\ja{[Should $\sgen$ in \eqref{eq:relationinverse} be $(\sgen)^{\star}$? Same for the text right after ~\eqref{eq:subMat}. ]}
It is worth noting that the value of $\JM$ computed via \eqref{eq:relationship}-\eqref{eq:J} depends on knowing the binding constraints $\setSgen$ and $\setSbra$ for given $(\f,\para,\sload)$.
We abuse notation slightly and let $\JM(\sload;\f,\para)$ be the Jacobian matrix when $(\f,\para,\sload)\in\setallr$ is known and let $\JM(\setSgen,\setSbra)$ be the Jacobian matrix when $(\setSgen,\setSbra)$ is known.
When it is clear from context or not relevant we simply use $\JM$.

\subsection{Range of OPF Derivative}
%In this subsection, we study the possible values of the Jacobian matrix $\JM(\sload;\f,\para)$ for $(\f,\para,\sload)\in\setallr$.
The previous subsection has shown that the value of $\JM(\sload;\f,\para)$ is equivalent to $\JM(\setSgen,\setSbra)$ for certain choice of $\setSgen$ and $\setSbra$.
The following theorem also implies the equivalence between the range  of $\JM(\sload;\f,\para)$ and $\JM(\setSgen,\setSbra)$.
\footnote{
Here, the range refers to the set of values that $\JM(\sload;\f,\para)$ or $\JM(\setSgen,\setSbra)$ could take,
rather than the column space of $\JM(\sload;\f,\para)$ or $\JM(\setSgen,\setSbra)$.
}
\begin{theorem}\label{thm:range-equiv}
\begin{align}\label{eq:range-equiv}
\nonumber
&\{ \JM(\sload;\f,\para)~|~(\f,\para,\sload)\in\setallr\}\\
\nonumber
 =&\{ \JM(\setSgen,\setSbra)~|~\setSgen\in\VertexG,\setSbra\in\Edge, \setSgen\perp\setSbra,\\
 &\quad\quad\quad\quad\quad~~ |\setSgen|+|\setSbra|=\nG-1\}.
\end{align}
\end{theorem}
Here, we use $\setSgen\perp\setSbra$ to denote that in \eqref{eq:opf1}, all the inequality constraints corresponding to $\setSgen$ and $\setSbra$, as well as equality constraints, are independent of each other.
Notice that the left hand side of \eqref{eq:range-equiv} is induced by the DC-OPF problem and hence involves physical parameters such as the cost function, generation and load. 
The right hand side, however, purely depends on the graph topology.
Theorem \ref{thm:range-equiv} shows the equivalence between the value ranges of $\JM(\sload)$ and $\JM(\setSgen,\setSbra)$.

We first provide  the following lemmas in order to build up to the final proof for Theorem \ref{thm:range-equiv}. We defer their proofs to Appendix \ref{app:proof_wc_reduction}.
\begin{lemma}\label{lm:find_f_para1}
For any $\setSgen\in\VertexG, \setSbra\in\Edge$ such that $|\setSgen|+|\setSbra|=\nG-1$ and $\setSgen\perp\setSbra$, there exist $(\fs,\paras,\sloads)\in\setall$
such that \eqref{eq:opf1} has unique solution and all the binding constraints at the solution point exactly correspond to $\setSgen$ and $\setSbra$.
\end{lemma}
\begin{IEEEproof}
See Appendix \ref{app:proof_wc_reduction}.
\end{IEEEproof}

\begin{lemma}\label{lm:find_f_para2}
For any $\setSgen\in\VertexG, \setSbra\in\Edge$ such that $|\setSgen|+|\setSbra|=\nG-1$ and $\setSgen\perp\setSbra$, there exist $\fss\in\setf$, $\parass\in\setparar(\fss)$ and an open ball $W\subseteq\setslr(\parass,\fss)$
such that all the binding constraints exactly correspond to $\setSgen$ and $\setSbra$ whenever $\sload\in W$.
\end{lemma}
\begin{IEEEproof}
See Appendix \ref{app:proof_wc_reduction}.
\end{IEEEproof}

Now we have all the ingredients for proving Theorem \ref{thm:range-equiv}.

\begin{IEEEproof}(Theorem \ref{thm:range-equiv})
%\ja{[Can we refer to corollary 3 in here?]} 
For any $(\f,\para,\sload)\in\setallr$, by definition the binding constraints $\setSgen$ and $\setSbra$ must satisfy $|\setSgen|+|\setSbra|=\nG-1$ and $\setSgen\perp\setSbra$. 
Thus the left hand side of \eqref{eq:range-equiv} is a subset of the right hand side of \eqref{eq:range-equiv}.
As for the opposite direction, Lemma \ref{lm:find_f_para2} implies for any $(\setSgen,\setSbra)$ such that $|\setSgen|+|\setSbra|=\nG-1$ and $\setSgen\perp\setSbra$ we can always find $(\f,\para,\sload)\in\setallr$ whose associated binding constraints exactly correspond to $(\setSgen,\setSbra)$. Hence the right hand side of \eqref{eq:range-equiv} is also a subset of the left hand side.
\end{IEEEproof}

The result of Theorem \ref{thm:range-equiv} also indicates there exists a surjection from $\setallr$ to the set $\{(\setSgen,\setSbra)~|~|\setSgen|+|\setSbra|=\nG-1, \setSgen\perp\setSbra\}$ and 
the derivative of the operator (depending on the parameters) and the Jacobian matrix (depending on the binding constraints combination) take the same value under such surjection.
If one is only interested in the range of $\JM(\sload;\f,\para)$ such as the worst-case analysis instead of the value at a specific point, then $\JM(\sload;\f,\para)$ and $\JM(\setSgen,\setSbra)$ may be used interchangeably. 
One benefit of studying $\JM(\setSgen,\setSbra)$ is it has a closed form expression and only depends on the graph topology of the system.
\section{Computation}\label{sec:comp}
%\ja{\textbf{Notation needs to be fixed. It's not currently consistent with the rest of the paper.}}
%Having established  sets on which the derivative of the $\OPF$ operator is guaranteed to exist, we now turn our attention to computation. 
In the previous section, we  provided a closed form expression for the Jacobian  $\JM=\partial_{\sload}\OPF(\sload)$ which  depends on the binding generators and branches. This expression  will be very useful in shedding light on  further properties of the sensitivity of the DC-OPF problem.
For instance, it helps us study the OPF sensitivity bounds in the ``worst case" \cite{AndZL20}, which provides privacy guarantees when releasing power flow data \cite{ZhoAL19}.
%On the other hand, the method above is restricted to the DC OPF problems and only reflects how the optimal solution changes when $\sload$ is perturbed. \ja{[What are we trying to say here?!]}
%\fz{[I simply removed the last sentence.]}

In this section, we will show how recent results on conic problem differentiation can be applied to the OPF operator, specifically in the case when one simply focuses on evaluating a derivative at a given operating point.
This method could provide the derivative of the optimal solution with respect to different system parameters, and could also be generalized to other power flow models.
For example, in Section \ref{subsec:ACOPF} we describe how these results can be applied to an AC OPF problem when a semidefinite relaxation of the power flow equations is considered. In such a setting we are unable to guarantee the existence of the derivative and we leave this to future work.

%Here we  provide a numerical method to compute derivatives of $\OPF$ with respect to arbitrary problem data in~\eqref{eq:opf1}.
%In Section \ref{subsec:general} we first rephrase a recent result \cite{ArgBBBM19} which is developed for a general conic program formulation. We then in Section \ref{subsec:DCOPF} and \ref{subsec:ACOPF} reformulate the DC and AC OPF problems into such conic program formulation so the theory can directly apply.

\subsection{Differentiating a General Conic Program}\label{subsec:general}
The method of computation we pursue largely follows that presented in~\cite{ArgBBBM19} which considers general convex conic optimization problems that are solved using the homogenous self-dual embedding framework~\cite{YeTM94, OdhCPB16}. Consider a standard  primal-dual pair written in conic form:
\begin{align*} &\text{(P)}~~
\begin{array}{lcl}
\underset{\cx,\cs}{\minimize}& &  {\cc}^\T \cx \\
\st && \cA\cx+\cs = \cb \\
&&  (\cx,\cs) \in \Real^n \times \cK,
\end{array}\\
&\text{(D)}~~
 \begin{array}{lcl}
  \underset{\cy,\ccr}{\minimize} &  & \cb^\T \cy \\
  \st&& \cA^\T \cy+\cc = \ccr\\ 
  & & (\ccr,\cy)\in \{0\}^n \times \cK^{\star}.
 \end{array}
\end{align*}
In this setting the problem data consists of the triple $(\cA,\cb,\cc)\in \Real^{m \times n} \times \Real^m \times \Real^n$. The primal variable is $\cx\in \Real^n$, the primal slack variable is $\cs \in \Real^m$, and the dual variable is $\cy\in \Real^m$, with $\ccr\in \Real^n$ the dual slack variable. The set $\cK$ in a non-empty, closed, convex cone with $\cK^{\star}$ its dual. Linear programming falls  into this class of conic problems by setting $\cK$ to be the positive orthant. 

The KKT conditions for primal-dual optimality are $\cA\cx+\cs=\cb$, $\cA^\T \cy +\cc= \ccr$, $\ccr=\zero$, $\cs \in \cK$, $\cy\in \cK^{\star}$, and $\cs^\T \cy=0$.  The homogenous self-dual embedding formulation is expressed as 
%\begin{subequations}\label{eq:selfdual}
\begin{align}
\text{find}  \quad &(\cu,\cv) \nonumber \\
\nonumber
\st \quad & \cv = \cQMat\cu \label{eq:sd_embed} \\
\quad & (\cu,\cv)\in \cC \times \cC^{\star} 
\end{align}
with cones $\cC = \Real^n \times \cK^{\star}\times \Real_+$ and its dual $\cC^{\star} = \{0\}^{n} \times \cK \times \Real_+$. The variables $\cu$ and $\cv$ correspond to variables in (P) and (D) and two augmented variables $\kappa$ and $\tau$, and satisfy the mapping:
\begin{align*}
&\underbrace{\left[\begin{array}{c} \ccr\\ \cs \\ \kappa \end{array} \right] }=
\underbrace{\left[\begin{array}{ccc} \zero & \cA^\T & \cc \\ -\cA & \zero & \cb \\ -\cc^\T & -\cb^\T & \zero \end{array} \right]}
\underbrace{\left[\begin{array}{c} \x\\ \y \\ \tau \end{array} \right]}, \quad (\tau, \kappa)\in \Real_+ \times \Real_+, \\
&\hspace{0.5cm}\cv  \hspace{2.25cm}\cQMat\hspace{1.9cm} \cu 
\end{align*}
which is exactly the affine constraint in~\eqref{eq:sd_embed}. 
{Using Minty's parametrization \cite{rockafellar2015convex}, we let $\cz\in \Real^{n+m+1}$  denote $\cu-\cv$, giving $\cu=\cP_{\cC}\cz$, and $\cv=-\cP_{-\cC^\star}\cz$.
Now reformulate \eqref{eq:sd_embed} in terms of $\cz$ as 
\begin{align}
\text{find}  \quad &\cz=(\cz_1\in\Real^n, \cz_2\in\Real^m, \cz_3\in\Real)\in\Real^{n+m+1} \nonumber \\
\nonumber
\st \quad & -\cP_{-\cC^\star}\cz = \cQMat\cP_{\cC}\cz \label{eq:sd_embed_z} \\
& \cz_3>0.
\end{align}

The \emph{solution map} is defined as $\cS: \Real^{m \times n}\times \Real^m \times \Real^{n}\rightarrow \Real^{2m+n}$ which ``pushes'' the problem data $(\cA,\cb,\cc)$ through optimization problem~\eqref{eq:sd_embed} to return $(\cx,\cy,\cs)$ -- the primal-dual solutions. As a functional,  we  can write $\cS = \psi \circ \phi \circ Q$. The function $Q$ constructs the skew-symmetric matrix {$\cQMat$} from $(\cA,\cb,\cc)$. The mapping  $\phi: \cQ \rightarrow \Real^{n+m+1}$ maps from the space of skew-symmetric matrices to solution {$\cz$} of the self-dual embedding \eqref{eq:sd_embed_z}.
Finally, $\psi: \Real^{n+m+1} \rightarrow \Real^n \times \Real^m \times \Real^m$ constructs the primal-dual solutions of (P) and (D) from the self-dual embedding solution, i.e. $(\cx,\cy,\cs) = \psi(\cz)$  where
\begin{equation*}
\psi(\cz) = (\cz_1, \cP_{\cK^{\star}}\cz_2, \cP_{\cK^{\star}}\cz_2-\cz_2 )/\cz_3
 \end{equation*}
with $\cz$  a solution of the self-dual embedding \eqref{eq:sd_embed_z}.
% and $\cP_{\cK^{\star}}$ is the projection operator onto the dual cone $\cK^{\star}$. For the positive orthant, given a vector $x$, the projector is $\cP_{\Real^n_+}x = \max\{x,0\}$
%applied in an element-wise manner.

The following result is taken from~\cite{ArgBBBM19}, it is essentially an application of the chain-rule and the implicit function theorem. Consider the perturbation in problem data, $(\dD \cA, \dD \cb, \dD \cc)$, and the derivative of the solution map, $\cD \cS/\cD(\cA,\cb,\cc)$, then the perturbation on the primal-dual solutions is evaluated from
\begin{align}
(\dD \cx, \dD \cy, \dD \cs) &= \frac{\cD \cS(\cA,\cb,\cc)}{\cD (\cA,\cb,\cc)}(\dD \cA, \dD \cb, \dD \cc) \\
 &= \frac{\cD \psi(\cz)}{\cD \cz}\frac{\cD \phi(\cQMat)}{\cD \cQMat} \frac{\cD Q(\cA,\cb,\cc)}{\cD (\cA,\cb,\cc)}(\dD \cA, \dD \cb, \dD \cc).
\end{align}
To evaluate the values of $(\dD \cx, \dD \cy, \dD \cs)$, we first derive the expression for $\dD \cz$ and then recover $(\dD \cx, \dD \cy, \dD \cs)$ from $\dD \cz$.
Numerically, \cite{ArgBBBM19} show that $\dD \cz = -\cM^{-1}\cg$, where
\begin{align}
\nonumber
\cM& = ((\cQMat-\eye)\DD \cP_{\cC}\cz+\eye)/\cz_3
\label{eq:M}\\
\nonumber
\cg &= \dD \cQMat \cP_{\cC} (\cz/\cz_3)\\
\dD \cQMat &= \left[\begin{array}{ccc} \zero & \dD \cA^\T &  \dD \cc \\ -\dD \cA & \zero & \dD \cb\\ -\dD \cc^\T & -\dD \cb^\T & \zero \end{array}\right].
\end{align}
Here we use $\DD$ instead of $\cD$ to denote the derivative of an operator when the arguments are clear from  context. 
Note that for large systems it may be preferable to not invert $\cM$ and instead solve a least squares problem. Finally, partition $\dD \cz$ conformally as $(\dD \cz_1, \dD \cz_2, \dD \cz_3)$ and compute
\begin{equation}\label{eq:dxdyds}
\left[ \begin{array}{c}\dD \cx \\  \dD \cy\\ \mathrm{d}\cs \end{array}   \right] = \left[ \begin{array}{c}
\mathrm{d}\cz_1 - (\mathrm{d}\cz_3)\cx \\
(\DD \cP_{\cK^{\ast}}(\cz_2))\dD \cz_2 - (\dD \cz_3)\cy\\
(\DD \cP_{\cK^{\ast}}(\cz_2))\dD \cz_2 - \dD \cz_2 - (\dD \cz_3)\cs
\end{array}
\right].
\end{equation}
The method outlined above provides us with more information than we have considered to this point. Specifically, it leverages information about the primal and dual conic forms and provides derivative information with respect to all problem data rather than just load changes.
%We have focussed on the case of sensitivity with respect to load changes, i.e. perturbations to the elements of the vector $\cb$. Thus, if we simply wish to numerically obtain the Jacobian, set $\dD \cA, \dD \cc$ and their transpositions to zero and $\dD \cb = \one$ in $\dD \cQMat$. It then follows that $\dD \cx = \cD{\sload}\OPF(\sload)$.

\subsection{DC Optimal Power Flow}\label{subsec:DCOPF}
The DC Optimal power flow problem~\eqref{eq:opf1} can easily be written in the form (P) by introducing the appropriate slack variables and taking $\cK =\{0\}^{\nGL+1}\times \Real_+^{2\nG+2\nE}$: 
%\begin{subequations}
\begin{eqnarray*}
\underset{\x:=[(\sgen)^\T,\ang^\T]^\T,\cs}{\text{minimize}}  && [\f^{\T}, \zero^{\T}]\x
\label{eq:opf-comp.a}\\
\text{ subject to }& & 
\left[
\begin{array}{c}
\Aeq\\
\Ain
\end{array}
\right]
 \x + \cs= 
 \left[
 \begin{array}{c}
 \beq\\
 \bin
 \end{array}
 \right]
\label{eq:opf-comp.b}\\
& & 
(\x,\cs)\in\Real^{\nG+\nGL}\times(\{0\}^{\nGL+1}\times \Real_+^{2\nG+2\nE})\nonumber\\
\label{eq:opf-comp.c}
\end{eqnarray*}
\label{eq:opf-comp}
%\end{subequations}
where $(\Aeq,\Ain,\beq,\bin)$ are as defined in \eqref{eq:standardAb}.
Noting that $\cK^{\star} = \Real^{\nGL+1}\times \Real_+^{2\nG+2\nE}$.

Here we note that the derivative  of the projection operator $\DD \cP_{\cC}$ appearing in~\eqref{eq:M} is decomposed as
\begin{equation*}
\DD \cP_{\Real^{\nG+\nGL}} \times \DD \cP_{\Real^{\nGL+1}} \times \DD \cP_{\Real^{2\nG+2\nE}_+} \times \DD \cP_{\Real_+}
\end{equation*}
and $\DD \cP_{\cK^*}$ appearing in \eqref{eq:dxdyds} is decomposed as
\begin{equation*}
\DD \cP_{\Real^{\nGL+1}} \times \DD \cP_{\Real^{2\nG+2\nE}_+} .
\end{equation*}
Specifically, $\DD \cP_{\Real_+}$ is differentiable everywhere but at $\{0\}$, elsewhere 
\begin{align*}
\DD \cP_{\Real_+}\x = \frac{1}{2}({\sign}(\x)+1).
\end{align*}

 \subsection{AC Optimal Power Flow}\label{subsec:ACOPF}
In this subsection we briefly outline how the methods described in the previous section extend seamlessly to a semidefinite programming based relaxation of the AC optimal power flow problem. Unlike with the DC case, we make no claim as to when the derivatives are guaranteed to exist.
%, or whether the perturbed solution allows us to reconstruct an optimal solution.

For AC OPF problems, the loads $\sload$ and generations $\sgen$ become complex numbers, where the real part denotes the real power and the imaginary part denotes the reactive power.
In \cite{Low2014convex}, the \emph{bus injection} AC OPF problem is formulated  as
\begin{subequations}
\begin{eqnarray}
\underset{\cW\in\Sspace_{\nGL}}{\text{minimize}}  &&~ \trace(\acC_0\cW)
\label{eq:opf-ac.a}\\
\text{ subject to }& & ~
\trace(\acC_i\cW)\leq \cb_i \quad i=1,2,\dots, m
\label{eq:opf-ac.b}\\
&&~\rank(\cW)=1
\label{eq:opf-ac.c}
\end{eqnarray}
\label{eq:opf-ac}
\end{subequations}
where $\Sspace_{\nGL}$ is the space of all the $\nGL\times\nGL$ positive semidefinite Hermitian matrices.
Matrices $\acC_i (i=0,1,\dots,m)$ are determined by the power system parameters such as admittances and the network topology.
The values $\cb_i (i=1,\dots,m)$ depend on both the load profile $\sload$ and system capacity limits, and $\cb$ is linear in $\sload$.
The optimal generation $\sgen=\OPF(\sload)$ is linear in the optimal solution $\cW^{\star}$ of \eqref{eq:opf-ac}.\footnote{Here we extend the notation $\OPF$ as the mapping that returns the optimal $\sgen$ for given $\sload$ based on the AC OPF problem.}
The task is  to now derive the derivative $\dD\sgen$ with respect to the perturbation $\sload$. Following the same arguments as the previous section, all that remains to be done is to numerically compute $\dD\cW$ for perturbations to $\cb$.

As \eqref{eq:opf-ac} is non-convex and thereby computationally challenging, the semidefinite relaxation is always applied by dropping the non-convex rank constraint \eqref{eq:opf-ac.c}.
For radial networks (i.e., when $\Graph$ is a tree), there are sufficient conditions under which the semidefinite relaxation is exact in the sense that it yields the same optimal solution as \eqref{eq:opf-ac}; see \cite{low2014convex2} and \cite{MolH19} for extensive references. % \cite{bose2015quadratically}. %The exactness of the relaxed solution and it's relationship with the topology of the network is an active area of research~\cite{MadSL14,SojL11}. 
The relaxed problem \eqref{eq:opf-ac.a}-\eqref{eq:opf-ac.b} is a semidefinite programming problem and thereby can be rewritten in the canonical form of a conic program as in (P) and (D).
The same technique in Section \ref{subsec:general} can be applied to numerically evaluate $\dD\cW$ -- the formulae for the derivative of the projection operator for the semidefinite cone can be found in~\cite{ArgBBBM19} and~\cite{MalS06}. It should be noted however, that perturbations to $\cW$ need not result in a rank-one solution.% Keep in for arXiv
\section{Illustrative Examples}\label{sec:example}
In this section, we use the IEEE 9-bus test network as an example to illustrate what the sets 
$(\setf,\setpara,\setparar,\setsl,\setslr)$
in Fig. \ref{Fig:Defs2} look like.
The topology of the network is shown in
Fig. \ref{Fig:network}. It has three generators (white circles) and 6 loads (black circles).
The susceptances (edge weights) of power lines are taken from the MATPOWER toolbox \cite{zimmerman2011matpower}. 
The system parameters are provided in Table \ref{tb:parameter}. 
The data for the capacity limits and the loads are either directly taken from MATPOWER or perturbed to satisfy our assumptions.

\begin{figure}
\centering
\includegraphics[width=0.9\columnwidth]{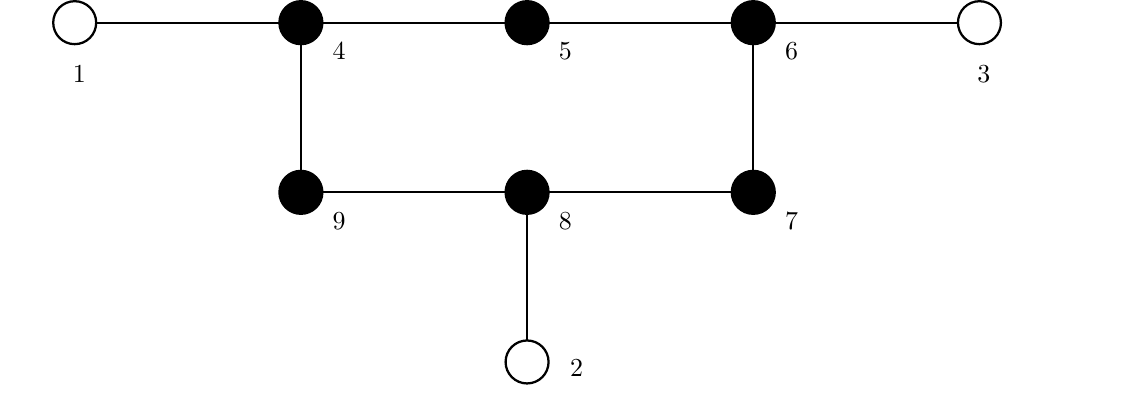}
\caption{Diagram of IEEE 9-bus test network. Generators are represented by white circles, while load buses are colored in black.}
\label{Fig:network}
\end{figure}

\begin{table*}[htp]
\caption{IEEE 9-bus Parameter Specification. The unit for the capacity limits and loads is 100 MW.}
\begin{center}
\begin{tabular}{|c|@{\hspace{0em}}c@{\hspace{0em}}|}
	\hline
	\rotatebox[origin=c]{90}{cost} &
	\begin{tabularx}{0.9\textwidth}{Y|Y|Y}
		$\f_1$ & $\f_2$ & $\f_3$\\
		\hline
		0.7191 & 0.5066 & 0.4758
	\end{tabularx}\\
	\hline
	\multirow{3}{*}{\rotatebox[origin=c]{90}{capacity limits}}&
	\begin{tabularx}{0.9\textwidth}{l|Y|Y|Y}
		$i\in\VertexG $ & 1 & 2 & 3 \\
		\hline
		$\genulim$ & 2.5679 & 3.0758 & 2.7743\\
		\hline
		$\genllim$ & 0.1392 & 0.1655 & 0.1171\\
		\hline
	\end{tabularx}\\
	&
	\begin{tabularx}{0.9\textwidth}{l|Y|Y|Y|Y|Y|Y|Y|Y}
		$\Edge$ & $(1,4)$ & $(4,5)$ & $(5,6)$ & $(3,6)$ & $(6, 7)$ & $(7, 8)$ & $(8, 9)$ & $(4, 9)$ \\
		\hline
		$\overline{\pflow}$ & 2.571 &  2.503 & 1.528 &  3.005 &  1.510 &  2.582 &  2.532 & 2.595\\
		\hline
		$\underline{\pflow}$ &    -2.503 & -2.544 & -1.538 & -3.077 & -1.580 & -2.519 & -2.545 & -2.565\\
		\hline
	\end{tabularx}\\
	& Visualization: the upper and lower bounds for branch $(2, 8)$\\
	\hline	
	\multirow{2}{*}{\rotatebox[origin=c]{90}{load}}&
	\begin{tabularx}{0.9\textwidth}{l|Y|Y|Y|Y}
		$i\in\VertexL $ & 5 & 6 & 8 & 9\\
		\hline
		$\sload$ & 0.90 & $10^{-10}$ & $10^{-10}$ & 1.25\\
		\hline
	\end{tabularx}\\
	& Visualization: the loads of buses $4$ and $7$\\
	\hline		
\end{tabular}
\end{center}
\label{tb:parameter}
\end{table*}%

First, we visualize and illustrate the sets $\Real_{+}^{\nG}$ and $\setf$ where the cost vector $\f$ resides. 
%It is noticed that the OPF problem always has the same optimal solutions when $\f$ is scaled by any strictly positive factor.
As we ignore the trivial case when $\f=\zero$, 
%to achieve the optimal visualization performance, 
we restrict $\f$ to the unit sphere for visual clarity.
As a result, $\Real_{+}^{\nG}$ is visualized by the blue region including the boundary and black curve segments shown in Fig. \ref{Fig:f}.
The black curve segments represent the set of $\f$ which may potentially make the OPF problem have multiple solutions or violate \eqref{eq:multiplier1}.
Thereby the blue region excluding the black curve segments is the restriction of a subset of $\setf$ onto the unit sphere. 
Figure \ref{Fig:f} provides a visualization that $\setf$ is dense in $\Real_{+}^{\nG}$ ($\nG=3$ in this example).
If the cost vector $\f$ is randomly chosen in $\Real_{+}^{\nG}$, then we will almost surely obtain a well-behaved $\f$ not aligned with the black curves.
In the rest of this example, we randomly pick $\f=[0.7191, 0.5066, 0.4758]^\T$, which is shown in the ``cost" sector in Table \ref{tb:parameter}, and visualized as the red point in Fig. \ref{Fig:f}.

\begin{figure}
\centering
\includegraphics[width=\columnwidth]{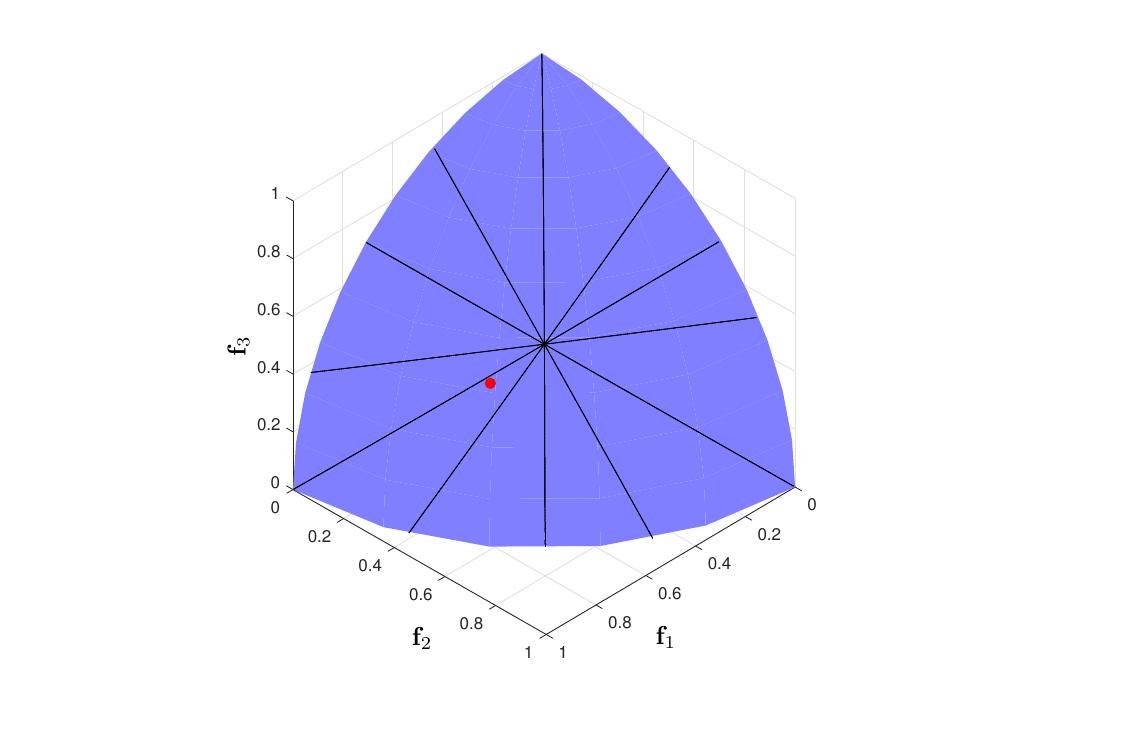}
\vspace{-2em}
\caption{Region for the cost vector $\f$. Black lines denote the values of $\f$ in $\Real_+^{\nG}\setminus\setf$. The red dot denotes the cost vector as indicated in Table \ref{tb:parameter} and used throughout this example.}
\label{Fig:f}
\end{figure}

\begin{figure}
\centering
\includegraphics[width=\columnwidth]{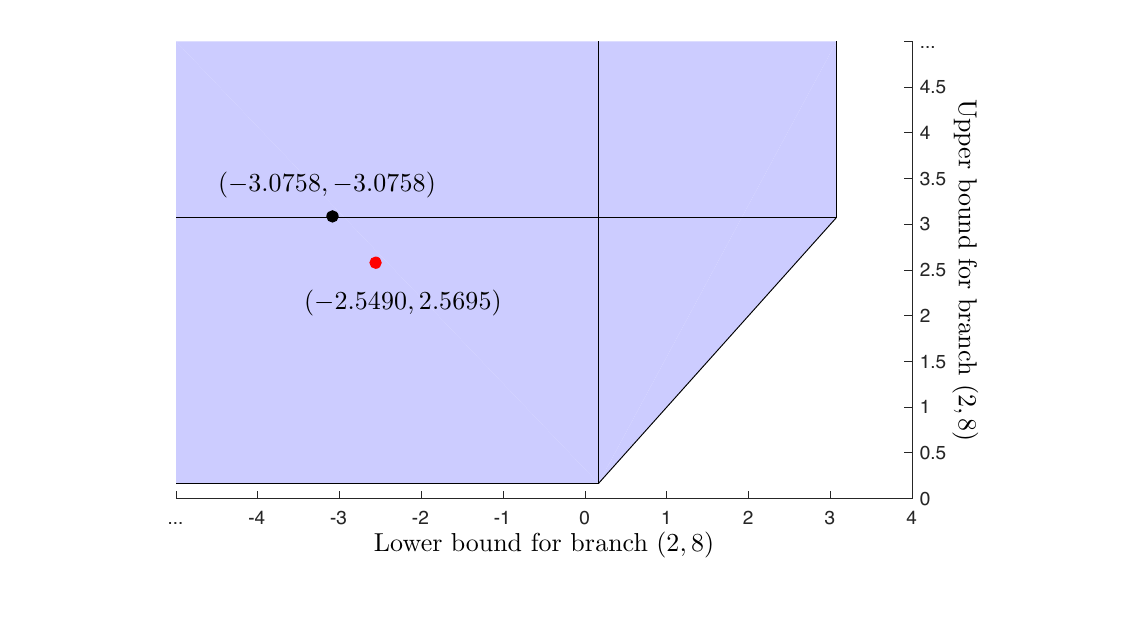}
\vspace{-2em}
\caption{Feasibility region for the power flow limits at branch $(2,8)$. The polytope (including its boundaries and the black lines) is (a slice of) $\setpara$. The set $\setparar(\f)$ is given by the purple region excluding the black lines and boundaries.  }
\label{Fig:para}
\end{figure}

We will now visualize the sets $\setpara$ and $\setparar(\f)$ for our choice of $\f$, and illustrate how different points in those two sets endow the OPF problem with different properties.
Consider that there are 3 generators and 9 branches in the network, and each generator and branch has both the upper and lower bounds for its generation and branch power flow,
the vector $\para$ thus has 24 dimensions. 
In order to make visualization possible, we fix all the capacity limits except for the power flow limits at branch $(2,8)$ as in Table \ref{tb:parameter}. 
A positive power flow at branch $(u,v)$ means that power is transmitted from $u$ to $v$.
Conversely, a negative value implies power is transmitted in the opposite direction.
Figure \ref{Fig:para} shows when $\f$ and other capacity limits are fixed, how the upper and lower bounds for branch $(2,8)$ affect the OPF operator.
In other words, Fig. \ref{Fig:para} visualizes a slice of sets $\setpara$ and $\setparar(\f)$.
The purple region, including the boundaries and black lines, is the slice of $\setpara$.
Picking any point in the purple region as the capacity limits for branch $(2,8)$, there exist some $\sload>0$ such that the constraints \eqref{eq:opf1.b}-\eqref{eq:opf1.e} are feasible.
However, for some points on the black lines or boundaries, the associated set $\setslr(\para,\f)$ might be not dense in $\setsl(\para)$.
We collect all the points in the purple region excluding the black lines and boundaries to form a slice of $\setparar(\f)$, which is dense in $\setpara$.
We now pick the red point in $\setparar(\f)$ (not on the black lines) and the black point in $\setpara\setminus\setparar(\f)$ (on the black line) as shown in Fig. \ref{Fig:para}, and will show their difference.
Recall that in Fig. \ref{Fig:Defs2}, we plot two points $\para_1\in\setparar(\f)$ and $\para_2\in\setpara\setminus\setparar(\f)$,
so the red point visualizes $\para_1$ while the black point visualizes $\para_2$.

\begin{figure*}
\centering
                  \begin{minipage}[b]{.7\columnwidth}
                    \centering
                    \includegraphics[width=1\columnwidth]{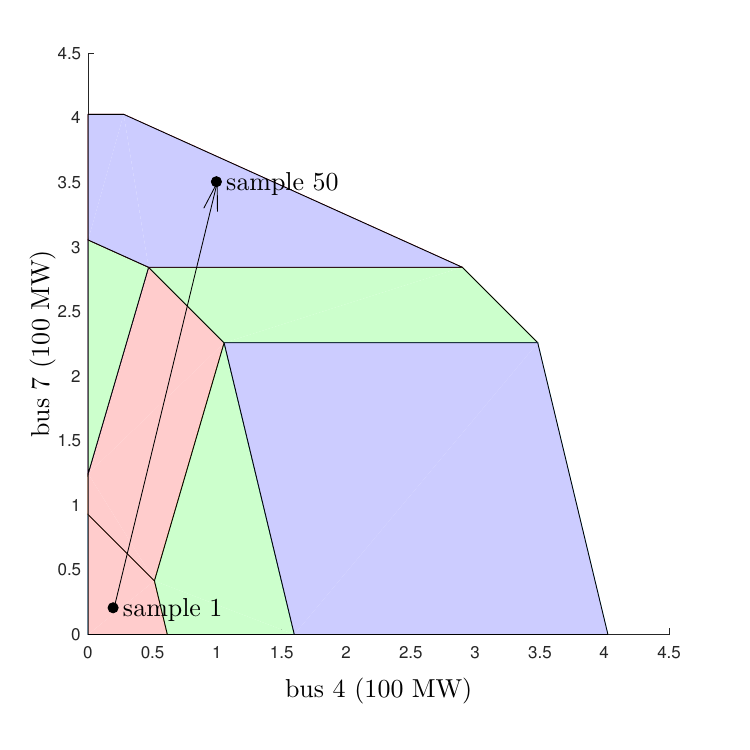}
                    \subcaption{When $\para$ is selected at the red point.}\label{Fig:load.a}
                  \end{minipage}%
                  \quad
                  \begin{minipage}[b]{.7\columnwidth}
                    \centering
                    \includegraphics[width=1\columnwidth]{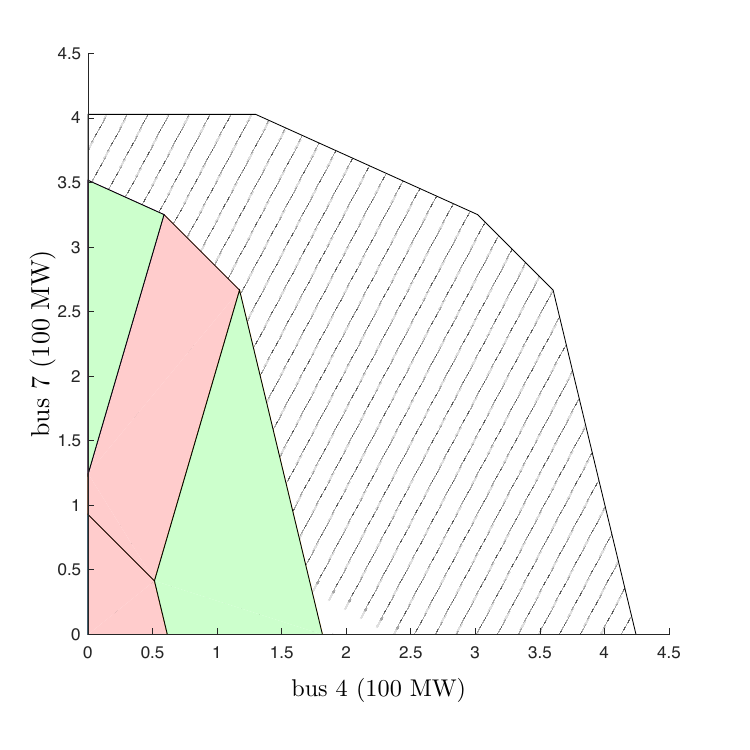}
                    \subcaption{When $\para$ is selected at the black point.}\label{Fig:load.b}
                  \end{minipage}
\caption{Region for the loads (bus 4 and bus 7).}
\label{Fig:load}
\end{figure*}

First, we pick the red point in Fig. \ref{Fig:para}, i.e., set the lower and upper bounds for branch $(2,8)$ at $(-2.5490,2.5695)$, respectively.
Since it is difficult to visualize all 6 loads, we fix buses $5$, $6$, $8$ and $9$ as in Table \ref{tb:parameter}, and visualize the region for buses $4$ and $7$ in Fig. \ref{Fig:load.a}.
The whole hexagon excluding the axes represents the slice of $\setsl(\para)$, within which any point corresponds to a load profile which makes the OPF problem feasible.
The whole region is further divided into seven colored subregions, and each of them refers to the set of load profiles under which the binding constraints of \eqref{eq:opf1} do not change.
In the interior of those subregions, there will be exactly $\nG-1=2$ independent binding inequality constraints.
Depending on the physical meaning of binding inequalities, we use three colors to distinguish different subregions.
Red indicates two binding constraints refer to two binding generators, green indicates one generator and one branch are binding, and purple indicates two binding branches.
Only the interior of those colored subregions contribute to the set $\setslr(\para,\f)$, which guarantees the number of and independence among all the binding constraints.
The operator $\OPF$ is also guaranteed to be differentiable when the loads are picked in $\setslr(\para,\f)$, and here the Jacobian matrix is given in Section \ref{subsec:Jacobian} in a closed form.
From Fig. \ref{Fig:load.a} we can see that when the red point is picked, the interior of all the subregions (i.e., $\setslr(\para,\f)$) is dense in the whole hexagon (i.e., $\setsl(\para)$).

Next, we pick the black point in Fig. \ref{Fig:para}, i.e., set the lower and upper bounds for branch $(2,8)$ at $(-3.0758,3.0758)$.
In this case, the whole hexagon contains a large chunk of shaded area. For the load profile in the shaded area, there might be more than $\nG-1=2$ binding inequality constraints, and all the binding constraints are not independent any more.
The Jacobian matrix we derived in Section \ref{subsec:Jacobian} is no longer valid.
As the shaded area is non-negligible, the interior of all the subregions is not dense in the whole hexagon any more.

Fortunately, both our theoretical proof and Fig. \ref{Fig:para} show that for almost all the capacity limits, they will behave like the red point in the above example and guarantee the independence among binding constraints for almost all the feasible load profiles.

\begin{figure*}
\centering
\includegraphics[width=1.35\columnwidth]{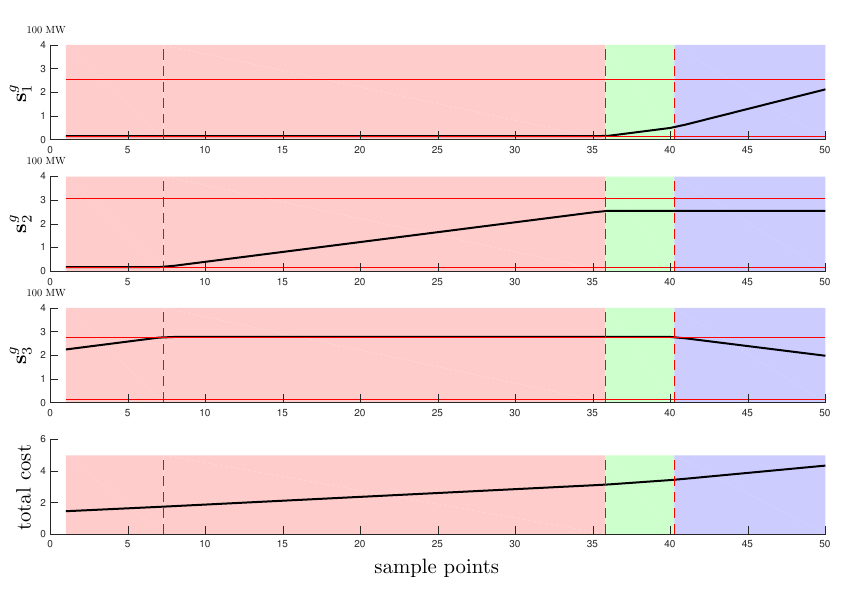}
\vspace{-1em}
\caption{Optimal generations and costs for different sample load profiles along the path plotted in Fig. \ref{Fig:load.a}. The solid red horizontal lines indicate the upper and lower bounds for the generation.}
\vspace{-2em}
\label{Fig:change}
\end{figure*}

Finally,  we consider a path in Fig. \ref{Fig:load.a} which goes through four different subregions, and pick 50 sample points along the path. 
Each sample point corresponds to a specific load profile for the power system.
In Fig. \ref{Fig:change}, we show how the optimal generations and costs change for those 50 sample load profiles.
In each subregion, the gradient of the optimal solution stays unchanged until the load profile enters a new subregion.
\section{Conclusion}
We presented an approach for analyzing a  linear program that solves the DC optimal power flow problem based on operator theoretic view of a linear program. Sets were defined upon which the OPF operator has a unique solution, is continuous, induce independent binding constraints, and the derivative exists (almost everywhere). Two equivalent perspectives on Jacobian matrix were given. The first was from the problem data and the second from knowledge of the binding constraints. A closed form expression of the Jacobian matrix is derived in terms of the binding constraints sets. Finally a numerical method based upon differentiating the solution map of a homogeneous self-dual conic program was described. 

It is hoped that this formulation will provide practitioners with new tools for analyzing the robustness of their networks. Simultaneously,  it opens up many interesting theoretical questions. In particular, studying AC optimal power flow problems from this perspective seems like a promising line of research. We are currently investigating how to compute the worst-case sensitivity of the DC-optimal power flow problem as this appears in a diverse range of applications including differential privacy, real-time optimal power flow problems, and locational marginal pricing.

\appendices
\section{Proof of Proposition \ref{fisdense}}\label{app:proof_fisdense}
We first define
\begin{subequations}
\begin{align}
\nonumber
\setf^{(1)}=\{&\f\in\Real_+^{\nG}~|~\forall \para\in\setpara, \forall \sload \in\setsl(\para), \eqref{eq:opf1}\text{~has}\\
&\text{unique optimal solution}\}\label{eq:omegaf.a}\\
\nonumber
\setf^{(2)}=\{&\f\in\Real_+^{\nG}~|~\forall \para\in\setpara, \forall \sload \in\setsl(\para), \text{all solutions}\\
&\text{of}~\eqref{eq:KKT}~\text{satisfy}~\eqref{eq:multiplier1}\}\label{eq:omegaf.b},
\end{align}\label{eq:omegaf}
\end{subequations}
then $\setf=\setf^{(1)}\cap\setf^{(2)}$.
For $\setS\subseteq\Edge$, $\setT\subseteq\VertexG$ such that $|\setS|+|\setT|\leq\nG-2$, we construct $\setQ(\setS,\setT)$ to be the set of $\f$ such that 
$\exists \taueq\in\Real^{\nGL+1},\mupm\in\Real^{\nE},\lambdapm\in\Real^{\nG}$ satisfying:
\begin{subequations}
\begin{eqnarray}
&& \bf{0}=M^{\T}\taueq+\CM\B\mupm
\label{eq:setQ.a}\\
&&  -\f=-[\taueq_1,\taueq_2,\cdots,\taueq_{\nG}]^{\T}+\lambdapm
\label{eq:setQ.b}\\
&&  \mupm_i\neq 0 \Rightarrow i\in\setS
\label{eq:setQ.c}\\
&&  \lambdapm_i\neq 0 \Rightarrow i\in\setT.
\label{eq:setQ.d}
\end{eqnarray}
\label{eq:setQ}
\end{subequations} 
When $\setS$ and $\setT$ are fixed, the vector $\CM\B\mupm$ takes value in an $|\setS|$ dimensional subspace. Since $\rank(\M)=\nGL$, the possible values of $\taueq$ must fall within an $|\setS|+1$ dimensional subspace.
Therefore, \eqref{eq:setQ.b} implies that $\f$ must be in an $|\setS|+1+|\setT|\leq\nG-1$ dimensional subspace, and hence $\inte(\clos(\setQ(\setS,\setT))=\emptyset$. 
Denote
\begin{align*}
\setQ_\cup:=
\Bigg(\bigcup\limits_{\tiny
\substack{\setS\subseteq\Edge,\setT\subseteq\VertexG\\
|\setS|+|\setT|\leq\nG-2 }}
\setQ(\setS,\setT)\Bigg),
\end{align*}
then, $\setQ_\cup\cap\Real_+^{\nG}$
is nowhere dense in $\Real_+^{\nG}$. 

On one hand, \eqref{eq:omegaf.b} and \eqref{eq:setQ} imply that 
\begin{align*}
\Real_+^{\nG}\backslash\setf^{(2)}&=
\{\f\in\Real_+^{\nG}~|&&\exists \para\in\setpara, \sload \in\setsl(\para), \text{one solution}\\
&&&\text{of}~\eqref{eq:KKT}~\text{violates}~\eqref{eq:multiplier1}\}\\
&\subseteq\setQ_\cup.
\end{align*}
Thereby, $\Real_+^{\nG}\backslash\setQ_\cup\subseteq\setf^{(2)}$.

On the other hand, we reformulate \eqref{eq:opf1} as
\begin{subequations}
\begin{eqnarray}
\underset{\x:=[(\sgen)^\T,\ang^\T]^\T}{\text{minimize}}  && [\f^{\T}, \zero^{\T}]\x
\label{eq:opf2.a}\\
\text{ subject to }& & 
\Aeq \x= \beq
\label{eq:opf2.b}\\
& & 
\Ain \x\leq \bin
\label{eq:opf2.c}
\end{eqnarray}
\label{eq:opf2}
\end{subequations}
where
\begin{subequations}\label{eq:standardAb}
\begin{eqnarray}
\nonumber
\Aeq:=\left[ \begin{array}{cc}
\zero^{1\times\nG} & \base_1 \\
\substack{-\eye^{\tiny\nG}\\\zero^{\tiny\nL\times\nG}} &
\CM\B\CM^\T
\end{array} \right], & 
\beq:=\left[ \begin{array}{c}
\zero^{(1+\nG)\times 1} \\ -\sload
\end{array} \right],\\
\end{eqnarray}
\begin{eqnarray}
\Ain:=\left[ \begin{array}{cc}
\zero^{\nE\times\nG} & \B\CM^\T\\
\zero^{\nE\times\nG} & -\B\CM^\T\\
{\eye}^{\nG} & \zero^{\nG\times\nGL} \\
 -{\eye}^{\nG} & \zero^{\nG\times\nGL}
\end{array} \right], & 
\bin:=\left[ \begin{array}{c}
\overline{\pflow} \\ -\underline{\pflow} \\ \genulim \\ -\genllim
\end{array} \right].
\end{eqnarray}
\end{subequations}
Geometrically, an LP has multiple optimal solutions if and only if 
the objective vector is normal to the hyperplane defined by equality constraints 
and the set of inequality constraints which are binding for all the optimal solutions (i.e., corresponding rows in $\Aeq$ and $\Ain$). 
We collect the rows in $\Ain$ which correspond to binding inequality constraints (for all the optimal solutions) and form a new matrix $\Aintilde$.
Formally, let $\setbinding$ be the set of indices $i$ such that the $i$\textsuperscript{th} row of $\Ain$ corresponds to a binding constraint for all the optimal solutions,
then $\Aintilde=\eye_{\setbinding}\Ain$.
In our case, the objective vector $[\f^{\T},\bf{0}^{\T}]^{\T}$ is an $\nG+\nGL$ dimensional vector,
thus the row space of $[\Aeq^\T,\Aintilde^\T]^\T$ must have dimension $\leq \nG+\nGL-1$ and $[\f^{\T},\bf{0}^{\T}]$ must be within this row space.
As $\Aeq$ has $\nGL+1$ linearly independent rows,
we can always find $\leq\nG-2$ independent rows of $\Aintilde$ to form a new matrix $\Aindbtilde$ such that
$[\Aeq^\T,\Aintilde^\T]^\T$ and $[\Aeq^\T,\Aindbtilde^\T]^\T$ share the same row space.
As a result,  $[\f^{\T},\bf{0}^{\T}]$ can be represented as the linear combination of rows in $[\Aeq^\T,\Aindbtilde^\T]^\T$,
and one can always find $(\setS,\setT,\taueq,\mupm,\lambdapm)$ satisfying \eqref{eq:setQ} and also $|\setS|+|\setT|\leq\nG-2$.
Hence $\Real_+^{\nG}\backslash\setf^{(1)}$ is also a subset of $\setQ_\cup$ and thus $\Real_+^{\nG}\backslash\setQ_\cup\subseteq\setf^{(1)}$.

Above all, 
$\Real_+^{\nG}\backslash\setQ_\cup\subseteq\setf^{(1)}\cap\setf^{(2)}=\setf$.
Since $\setQ_\cup\cap\Real_+^{\nG}$ is nowhere dense in $\Real_+^{\nG}$,
$\setf$ is dense in $\Real_+^{\nG}$.

\section{preliminaries for the Proof of Theorem \ref{limitsdense}}\label{app:proof_setparaissolid}
The following results are used in the proof of Theorem~\ref{limitsdense}. Together they show that any subset of $\setpara$ which can be covered by a finite union of subspaces of lower dimensions must be nowhere dense in $\setpara$. Further, it shows that  ``good" $\para$ always form a dense subset of $\setpara$.
\begin{proposition}\label{setparaissolid}
The set $\setpara$ satisfies
$\clos(\inte(\setpara)) = \clos(\setpara)$.
\end{proposition}
%\begin{IEEEproof}
%See Appendix \ref{app:proof_setparaissolid}.
%\end{IEEEproof}
%\fz{This proposition, together with Lemma \ref{lm:dense} in Appendix \ref{app:proof_limitsdense}, will imply that 
%any subset of $\setpara$ which can be covered by a finite union of subspaces of lower dimensions must be nowhere dense in $\setpara$.}
%Further, it can show ``good" $\para$ always form a dense subset of $\setpara$.
\begin{IEEEproof}
As it is trivial that $\clos(\inte(\setpara))\subseteq \clos(\setpara)$, we only need to show $\clos(\setpara)\subseteq \clos(\inte(\setpara))$.
It is sufficient to show that $\forall \para\in\setpara$,  there exists a sequence $(\para_{(n)})_{n=1}^{\infty}$ such that $\lim\limits_{n\to\infty}\para_{(n)}=\para$ and for each $n\in\Zset_+$ 
there is an open neighborhood $U(\para_{(n)})\ni\para_{(n)}$ such that $U(\para_{(n)})\subseteq\setpara$. 
The reason we only need to show this, is that, by definition, any point in $\clos(\setpara)$ is the limit of a sequence of points in $\setpara$. If any point in $\setpara$ is further the limit of a sequence of points in $\inte(\setpara)$, 
then any point in $\clos(\setpara)$ can also be represented as the limit of a sequence of points in $\inte(\setpara)$. That is to say, $\clos(\setpara)\subseteq \clos(\inte(\setpara))$. Next we prove $\forall \para\in\setpara$,  such sequence $(\para_{(n)})_{n=1}^{\infty}$ exists.

%\ja{\textbf{[We had this discussion in NETLAB a few weeks ago. If there is a standard analysis reference we can give that uses proofs of this type then could we provide it here?]}}

First, we observe that \eqref{eq:opf1.b}-\eqref{eq:opf1.e} implies the branch power flow $\pflow:=\B\CM^\T\ang$ satisfies 
\begin{align*}
\pflow=\B\CM^\T(\CM\B\CM^\T)^{\dagger}
\left[
\begin{array}{c}
\sgen\\
-\sload
\end{array}
\right].
\end{align*}
We use $\rho$ to denote the matrix norm of $\B\CM^\T(\CM\B\CM^\T)^{\dagger}$ induced by the $\ell_1$ vector norm:
\begin{align*}
\rho:=\underset{\x\neq\zero}{\rm maximize}\frac{\|\B\CM^\T(\CM\B\CM^\T)^{\dagger}\x\|_1}{\|\x\|_1}.
\end{align*}
Now consider any $\hat{\para}=[(\genulim)^{\T},(\genllim)^{\T},\overline{\pflow}^{\T},\underline{\pflow}^{\T}]^{\T}\in\setpara$ with $\genllim\geq 0$, 
and there exists $(\sgenh,\sloadh)$ such that $\sloadh> 0$ and \eqref{eq:opf1.b}-\eqref{eq:opf1.e} are satisfied (with associated branch power flow $\hat{\pflow}$).
Then we construct $\para_{(n)}$ as
\begin{align*}
\para_{(n)}=\left[
\begin{array}{c}
\genulim+\frac{3}{n}\one_{\nG}\\
\genllim+\frac{1}{n}\one_{\nG}\\
\overline{\pflow}+5\rho\frac{\nG}{n}\one_{\nE}\\
\underline{\pflow}-5\rho\frac{\nG}{n}\one_{\nE}
\end{array}
\right]
\end{align*}
and its open neighborhood
\begin{align*}
U(\para_{(n)})=\left\{\para\middle|
-\left[
\begin{array}{c}
\frac{1}{n}\one_{2\nG}\\
\rho\frac{\nG}{n}\one_{2\nE}
\end{array}
\right]
<\para-\para_{(n)}<
\left[
\begin{array}{c}
\frac{1}{n}\one_{2\nG}\\
\rho\frac{\nG}{n}\one_{2\nE}
\end{array}
\right]
\right\}.
\end{align*}
Clearly, $\para_{(n)}$ converges to $\hat{\para}$ as $n\to\infty$. Next, we are going to prove that for any $\para\in U(\para_{(n)})$, we have $\para\in\setpara$.
For the convenience of notation, we use $\genulim(\para)$, $\genllim(\para)$, $\overline{\pflow}(\para)$, $\underline{\pflow}(\para)$ to denote the corresponding part in $\para$.
Since for any $i\in\VertexG$
\begin{align*}
\genllim_i(\para)&>\genllim_i(\para_{(n)})-\frac{1}{n}
=\genllim_i(\hat{\para})+\frac{1}{n}-\frac{1}{n}
=\genllim_i(\hat{\para})\geq 0,
\end{align*}
we only need to check if there exist $(\sgen,\sload)$ such that \eqref{eq:opf1.b}-\eqref{eq:opf1.e}  are satisfied and $\sload>0$.
We construct $\sgen=\sgenh+\frac{2}{n}\one_{\nG}$ and $\sload=\sloadh+\frac{2\nG}{n\nL}\one_{\nL}$,
then it is clear that $\sload\geq\sloadh> 0$.
Since $\one_{\nG}^\T\sgen=\one_{\nG}^\T\sgenh+\frac{2\nG}{n}=\one_{\nL}^\T\sloadh+\frac{2\nG}{n}=\one_{\nL}^\T\sload$, the constructed generation and load are balanced so 
\eqref{eq:opf1.c} is satisfied for some $\ang$. Further, we can always shift $\ang$ to make $\ang_1=0$ and \eqref{eq:opf1.b} is thereby satisfied. Next, we can check 
\begin{align*}
\sgen&=\sgenh+\frac{1}{n}\one_{\nG}+\frac{1}{n}\one_{\nG}
\geq\genllim(\para_{(n)})+\frac{1}{n}\one_{\nG}>\genllim(\para)\\
\sgen&=\sgenh+\frac{3}{n}\one_{\nG}-\frac{1}{n}\one_{\nG}
\leq\genulim(\para_{(n)})-\frac{1}{n}\one_{\nG}<\genulim(\para),
\end{align*}
thus \eqref{eq:opf1.d} is satisfied. Finally, 
\begin{align*}
\pflow=&\B\CM^\T(\CM\B\CM^\T)^{\dagger}
\left[
\begin{array}{c}
\sgenh+\frac{2}{n}\one_{\nG}\\
-\sloadh-\frac{2\nG}{n\nL}\one_{\nL}
\end{array}
\right]\\
=&\hat{\pflow}+\B\CM^\T(\CM\B\CM^\T)^{\dagger}
\left[
\begin{array}{c}
\frac{2}{n}\one_{\nG}\\
-\frac{2\nG}{n\nL}\one_{\nL}
\end{array}
\right],
\end{align*}
so 
\begin{align*}
\pflow_i&\geq\hat{\pflow}_i-\left\|
\B\CM^\T(\CM\B\CM^\T)^{\dagger}
\left[
\begin{array}{c}
\frac{2}{n}\one_{\nG}\\
-\frac{2\nG}{n\nL}\one_{\nL}
\end{array}
\right]
\right\|_1\\
&\geq\underline{\pflow}_i(\hat{\para})-\rho\cdot\left(\frac{2}{n}\nG+\frac{2\nG}{n\nL}\nL\right)\\
&=\underline{\pflow}_i(\hat{\para})-4\rho\frac{\nG}{n}
=\underline{\pflow}_i(\para_{(n)})+\rho\frac{\nG}{n}>\underline{\pflow}_i(\para)
\end{align*}
and similarly $\pflow_i<\overline{\pflow}_i(\para)$. As the result, we have $\underline{\pflow}(\para)\leq\pflow\leq\overline{\pflow}(\para)$
and \eqref{eq:opf1.e} is satisfied.

Above all, we have shown that there exist $(\sgen,\sload)$ such that \eqref{eq:opf1.b}-\eqref{eq:opf1.e}  are satisfied and $\sload>0$. Thus $\para\in\setpara$ and there must be $U(\para_{(n)})\subseteq\setpara$.
\end{IEEEproof}
\begin{lemma}\label{lm:dense}
Suppose the set $\setS\subseteq\Real^{n}$ satisfies the condition that $\clos(\inte(\setS))=\clos(\setS)$, and $\setT$ is an affine hyperplane with dimension strictly less than $n$.
Then $\setT$ is nowhere dense in $\setS$.
\end{lemma}
\begin{IEEEproof}
If not, then by definition, in the relative topology of $\setS$, we have $\inte(\setT)\neq\emptyset$ since $\setT$ is closed. 
Pick any point $\x\in\inte(\setT)$, there must be an $n$-dimensional open ball $U$ with radius $r$ centered at $\x$ such that $\x\in U\cap\setS\subseteq\setT$.
In the $n$-dimensional Euclidean topology, since $\clos(\inte(\setS))=\clos(\setS)$, 
there must be a point $\x_1\in\setS$ such that $|\x-\x_1|\leq r/2$ and there is an $n$-dimensional open ball $U_1$ centered at $\x_1$ and have radius $< r/2$ satisfying $U_1\subseteq\setS$.
Clearly, $U_1\subseteq U$ as well, and thereby $U_1\subseteq U\cap\setS\subseteq \setT$. However, $\setT$ is an affine hyperplane with dimension strictly less than $n$, and there is the contradiction.
\end{IEEEproof}

\section{Proof of Theorem \ref{limitsdense}}\label{app:proof_limitsdense}
Our strategy is to construct the set $\setparar(\f)$ first, then prove $\clos(\inte(\setsl(\para)))=\clos(\setsl(\para))$, and finally show that $\setslr(\para,\f)$  is dense in $\setsl(\para)$.

Consider the power flow equations below:
\begin{align}\label{eq:constraints}
\mathbf{T}\ang:=
\left[
\begin{array}{c}
\CM\B\CM^{\T}\\
\B\CM^{\T}
\end{array} 
\right]\ang=
\left[
\begin{array}{c}
\sgen\\
-\sload\\
\pf
\end{array} 
\right].
\end{align}
%Here, $\pf\in\Real^{\nE}$ is the power flow for branches.
Proposition \ref{fisdense} and Assumption \ref{A:vectorf} show that there will always be at least $\nG-1$ binding inequality constraints as each non-zero multiplier will force one inequality constraint to be binding.
A constraint being binding means some $\sgen_i$ equals either $\genulim_i$ or $\genllim_i$ (as in the upper $\nG$ rows in \eqref{eq:constraints}),
or some $\pf_i$ equals either $\overline{\pflow}_i$ or $\underline{\pflow}_i$ (as in the lower $\nE$ rows in \eqref{eq:constraints}).
We have $\rank(\mathbf{T})=\nGL-1$.
%For $\setS\subseteq\VertexG\cup[\nGL+1,\nGL+\nE]$, let $\mathbf{T}_{\setS}$ be the matrix of $\mathbf{T}$'s rows whose indexes are in $\setS$, and the rows are ordered by the indexes.
%Similarly, let $\mathbf{T}_{\setS}$ be the matrix of $\mathbf{T}$'s rows whose indexes are in $\setS$.
We use the following procedure to construct the set $\setparar$.
\begin{enumerate}[I.]
\item $\setparar \leftarrow \setpara \backslash\big(\big(\bigcup\limits_{i\in\VertexG}\{\para~|~\para_i=\para_{i+\nG}\}\big)\cup\big(\bigcup\limits_{i-2\nG\in\Edge}\{\para~|~\para_i=\para_{i+\nE}\}\big)\big)$
\item For each $\setS\subseteq\VertexG\cup[\nGL+1,\nGL+\nE]$, construct $\mathbf{T}_{\setS}$.
\begin{enumerate}[a)]
\item If $\rank(\mathbf{T}_{\setS})=|\setS|$, then continue to another $\setS$.
\item If $\rank(\mathbf{T}_{\setS})<|\setS|$, then consider
\begin{align}\label{eq:Gamma}
\Gamma:=\prod\limits_{\mathclap{\tiny \substack{i\in\setS\cap\VertexG}}}~
\{\base_{\scriptscriptstyle i},\base_{\scriptscriptstyle\nG+i}\}\times
\prod\limits_{\mathclap{\tiny
\substack{j\in\Edge\\j+\nGL\in\setS}}}~
\{\base_{\scriptscriptstyle 2\nG+j},\base_{\scriptscriptstyle 2\nG+\nE+j}\}.
\end{align}
%where $\base_{\scriptscriptstyle m}$ denotes the standard base for the $m$\textsuperscript{th} coordinate.
%In the context of \eqref{eq:Gamma}, $\base_{\scriptscriptstyle m}\in\Real^{2\nG+2\nE}$.
Now update $\setparar$ as
\begin{align}\label{eq:setupdate}
\setparar \leftarrow\setparar \setminus
\bigcup\limits_{\gamma\in\Gamma}\big\{\para~|~\exists\ang,~\suchthat~\gamma^{\T}\para=\mathbf{T}_{\setS}\ang\big\}.
\end{align}
\end{enumerate}
\item Return $\setparar$.
\end{enumerate}

\begin{figure*}
\begin{align}\label{eq:Gamma}
\nonumber
\Gamma=\{&\base_{\scaleto{1}{4pt}},\base_{\scaleto{\nG+1}{4pt}}\}\times\{\base_{\scaleto{4}{4pt}},\base_{\scaleto{\nG+4}{4pt}}\}\times\{\base_{\scaleto{2\nG+2}{4pt}},\base_{\scaleto{2\nG+\nE+2}{4pt}}\}\\
\nonumber
=\{&(\base_{\scaleto{1}{4pt}}, \base_{\scaleto{4}{4pt}}, \base_{\scaleto{2\nG+2}{4pt}}),
(\base_{\scaleto{1}{4pt}}, \base_{\scaleto{4}{4pt}}, \base_{\scaleto{2\nG+\nE+2}{4pt}}),
(\base_{\scaleto{1}{4pt}}, \base_{\scaleto{\nG+4}{4pt}}, \base_{\scaleto{2\nG+2}{4pt}}),
(\base_{\scaleto{1}{4pt}}, \base_{\scaleto{\nG+4}{4pt}}, \base_{\scaleto{2\nG+\nE+2}{4pt}}),\\
\nonumber
&(\base_{\scaleto{\nG+1}{4pt}}, \base_{\scaleto{4}{4pt}}, \base_{\scaleto{2\nG+2}{4pt}}),
(\base_{\scaleto{\nG+1}{4pt}}, \base_{\scaleto{4}{4pt}}, \base_{\scaleto{2\nG+\nE+2}{4pt}}),
(\base_{\scaleto{\nG+1}{4pt}}, \base_{\scaleto{\nG+4}{4pt}}, \base_{\scaleto{2\nG+2}{4pt}}),
(\base_{\scaleto{\nG+1}{4pt}}, \base_{\scaleto{\nG+4}{4pt}}, \base_{\scaleto{2\nG+\nE+2}{4pt}})
\}\\
\nonumber
=\{
&[\base_{\scaleto{1}{4pt}}~ \base_{\scaleto{4}{4pt}}~ \base_{\scaleto{2\nG+2}{4pt}}],
[\base_{\scaleto{1}{4pt}}~ \base_{\scaleto{4}{4pt}}~ \base_{\scaleto{2\nG+\nE+2}{4pt}}],
[\base_{\scaleto{1}{4pt}}~ \base_{\scaleto{\nG+4}{4pt}}~ \base_{\scaleto{2\nG+2}{4pt}}],
[\base_{\scaleto{1}{4pt}}~ \base_{\scaleto{\nG+4}{4pt}}~ \base_{\scaleto{2\nG+\nE+2}{4pt}}],\\
&[\base_{\scaleto{\nG+1}{4pt}}~ \base_{\scaleto{4}{4pt}}~ \base_{\scaleto{2\nG+2}{4pt}}],
[\base_{\scaleto{\nG+1}{4pt}}~ \base_{\scaleto{4}{4pt}}~ \base_{\scaleto{2\nG+\nE+2}{4pt}}],
[\base_{\scaleto{\nG+1}{4pt}}~ \base_{\scaleto{\nG+4}{4pt}}~ \base_{\scaleto{2\nG+2}{4pt}}],
[\base_{\scaleto{\nG+1}{4pt}}~ \base_{\scaleto{\nG+4}{4pt}}~ \base_{\scaleto{2\nG+\nE+2}{4pt}}]
\}
\end{align}
\begin{spacing}{0.0}
\hrulefill
\end{spacing}
\end{figure*}

In the above procedure, an $n$-tuple of vectors is also regarded as a matrix of $n$ columns and the product in \eqref{eq:Gamma} is Cartesian product.
\footnote{Hence, each $\gamma\in\Gamma$ can also be regarded as a $(2\nG+2\nE)$-by-$|\setS|$ matrix.
For instance, if we have $\setS=\{1,4,\nGL+2\}$, then \eqref{eq:Gamma} shows a set of 8 elements and each element is a $(2\nG+2\nE)$-by-$3$ matrix.}
Since $\gamma\in\Gamma$ is of rank $|\setS|$ and $\mathbf{T}_{\setS}\ang$ with $\ang\in\Real^{\nGL}$ defines a subspace of $\leq|\setS|-1$ dimensions,
each set of $\{\para~|~\exists\ang,~\suchthat~\gamma^{\T}\para=\mathbf{T}_{\setS}\ang\}$ in \eqref{eq:setupdate} is a subspace with dimension strictly lower than $2\nG+2\nE$, 
and is thereby nowhere dense in $\setpara$ by Lemma \ref{lm:dense}.
Similarly, the sets $\{\para~|~\para_i=\para_{i+\nG}\}$ for $i\in\VertexG$ and $\{\para~|~\para_i=\para_{i+\nE}\}$ for $i-2\nG\in\Edge$ are also nowhere dense.
As a result, we have that $\setparar$ is dense in $\setpara$.
It is sufficient to show that two conditions in Proposition \ref{limitsdense} are satisfied.

To show $\clos(\inte(\setsl(\para)))=\clos(\setsl(\para))$, it is sufficient to prove that fix $\para\in\setparar$, $\forall\sloadh\in\setsl(\para)$, 
there exists a sequence $(\sload_{(n)})_{n=1}^{\infty}$ such that $\lim\limits_{n\to\infty}\sload_{(n)}=\sloadh$ and each $\sload_{(n)}$ has an open neighborhood $U(\sload_{(n)})$ such that $U(\sload_{(n)})\subseteq\setsl(\para)$. 
By definition, there exists $\sgenh$ and $\angh$ such that \eqref{eq:opf1.b}-\eqref{eq:opf1.e} are satisfied for $\sloadh$. We also use $\pflowh$ to denote the branch power flow associated with $(\sgenh,\sloadh)$.
Here we overload $\setS\subseteq\VertexG\cup[\nGL+1,\nGL+\nE]$ to denote the indices of all the binding inequality constraints for $(\sgenh,\sloadh)$. 
\footnote{In this section, the index of a constraint associated with generator $i$ (either the upper or lower bounds) is $i$ and 
the index of a constraint associated with branch $i$ (either the upper or lower bounds) is $i+\nGL$. Step I in the procedure constructing $\setparar$ guarantees that a generator or branch cannot reach the upper and lower bound at the same time.}
By construction, we have $\rank(\mathbf{T}_{\setS})=|\setS|\leq\rank(\mathbf{T})=\nGL-1$.
There are two situations to discuss: $|\setS|=0$ and $1\leq|\setS|\leq\nGL-1$.

In the first case, if $|\setS|=0$, then let $\rho_1$ be the matrix norm of 
\begin{align*}
\mathbf{T}_{\VertexG\cup[\nGL+1,\nGL+\nE]}
\left[
\begin{array}{c}
\base_1^\T\\
\mathbf{T}_{\VertexL}
\end{array}
\right]^{\dagger}
\end{align*}
induced by the $\ell_1$ vector norm. Let
\begin{align}
\nonumber
\epsilon_1&=
\min\Bigg\{
\mind\limits_{\tiny\substack{i\in\VertexG\\\sgenh_i>\genllim_i}} \sgenh_i-\genllim_i,
\mind\limits_{\tiny\substack{i\in\VertexG\\\genulim_i>\sgenh_i}} \genulim_i-\sgenh_i,\\
&\quad\quad\quad\quad\mind\limits_{\tiny\substack{i\in\Edge\\\pflowh_i>\underline{\pflow}_i}} \pflowh_i-\underline{\pflow}_i,
\mind\limits_{\tiny\substack{i\in\Edge\\\overline{\pflow}_i>\pflowh_i}} \overline{\pflow}_i-\pflowh_i
\Bigg\}\label{eq:eps1}\\
%\end{align}
%and %\ja{\textbf{[This notation is confusing because we're using min as minimum element of a set and min as short hand for minimize. Can we use min. for minimize and then add an explanation here?]}}
%\begin{align}\label{eq:eps2}
\epsilon_2 &= \mind\limits_{i}\sloadh_i.\label{eq:eps2}
\end{align}
Here, we have used $\mind$ as short hand for $\text{minimize}$.
Now we can construct $\sload_{(n)}\equiv\sloadh$, and 
\begin{align*}
U(\sload_{(n)})=\left\{
\sload\middle|
|\sload-\sload_{(n)}|<\frac{1}{2}\min\{\frac{\epsilon_1}{\nL\rho_1},\epsilon_2\}\one_{\nL}
\right\}.
\end{align*}
It is trivial that $\lim\limits_{n\to\infty}\sload_{(n)}=\sloadh$. For any $\sload\in U(\sload_{(n)})$, we have
\begin{align*}
\sload>\sload_{(n)}-\frac{1}{2}\epsilon_2\one_{\nL}=\sloadh-\Big(\frac{1}{2} \mind\limits_{i}\sloadh_i\Big)\one_{\nL}>0.
\end{align*}
Further, we will show that for
\begin{align}
\ang=\left[
\begin{array}{c}
\base_1^\T\\
\mathbf{T}_{\VertexL}
\end{array}
\right]^{\dagger}
\left[
\begin{array}{c}
0\\
-\sload
\end{array}
\right], 
\sgen=\mathbf{T}_{\VertexG}\ang,
\label{eq:construct_sg_ang}
\end{align}
\eqref{eq:opf1.b}-\eqref{eq:opf1.e} are satisfied. 
%In \eqref{eq:construct_sg_ang}, $\setS':=\{g(i)|i\in\setS\}$ where
%\begin{align*}
%g(i):=\left\{\small
%\begin{array}{ll}
%i,~&i\in\VertexG,~\sgenh_i=\genulim_i\\
%i+\nG,&i\in\VertexG,~\sgenh_i=\genllim_i\\
%i-\nGL+2\nG,&i-\nGL\in\Edge,~\pflowh_{i-\nGL}=\overline{\pflow}_{i-\nGL}\\
%i-\nGL+2\nG+\nE,&i-\nGL\in\Edge,~\pflowh_{i-\nGL}=\underline{\pflow}_{i-\nGL}
%\end{array}
%\right..
%\end{align*}
Clearly \eqref{eq:construct_sg_ang} implies $\base_1^\T\ang=0$, $\mathbf{T}_{\VertexL}\ang=-\sload$ and $\mathbf{T}_{\VertexG}\ang=\sgen$, which are equivalent to \eqref{eq:opf1.b}, \eqref{eq:opf1.c}.
For \eqref{eq:opf1.d}, as no $\sgenh_i$ reaches any bound, we have 
\begin{align*}
|\sgen_i-\sgenh_i|&\leq\rho_1\|\sload-\sloadh\|_1
\leq\rho_1\cdot\frac{1}{2}\nL\frac{\epsilon_1}{\nL\rho_1}=\frac{1}{2}\epsilon_1,
\end{align*}
and thereby $\sgen_i$ is still strictly between the bounds and stays feasible.
Similarly, the branch flow $\pflow$ is also within the upper and lower bounds, and \eqref{eq:opf1.e} is also satisfied.
As a result, $\sload\in\setsl(\para)$ and thus $U(\sloadh)\subseteq\setsl(\para)$.

In the second case, we have $1\leq|\setS|\leq\nGL-1$, then define
\begin{align*}
\mathbf{T}'(\setR):=\left[
\begin{array}{c}
\base_1^\T\\
\mathbf{T}_{\setS}\\
\mathbf{T}_{\setR}
\end{array}
\right],
~\text{for}~
\setR\subseteq\VertexL.
\end{align*}
Let $\setR^*=\argmin\limits_{\setR:
\rank(\mathbf{T}'(\setR))=\rank(\mathbf{T}'(\VertexL))} |\setR|$.
If there are multiple  $\setR$ that minimize $|\setR|$ then pick any one of them.
There are two simple observations:
\begin{itemize}
\item All rows of matrix $\mathbf{T}'(\setR^*)$ are independent.
\item All rows of $\mathbf{T}_{\VertexL}$ are in the row space of $\mathbf{T}'(\setR^*)$.
\end{itemize}

We further define
\begin{align*}
\mathbf{T}''(\setT):=\left[
\begin{array}{c}
\mathbf{T}'(\setR^*)\\
\mathbf{T}_{\setT}
\end{array}
\right],
\end{align*}
for $\setT\subseteq(\VertexG\cap[\nGL+1,\nGL+\nE])\backslash\setS$.
Let 
\begin{align*}
\setT^*=\argmin\limits_{\setT:
\rank(\mathbf{T}''(\setT))=\rank(\mathbf{T}''((\VertexG\cap[\nGL+1,\nGL+\nE])\backslash\setS))} |\setT|.
\end{align*}
Likewise, if there are multiple such $\setT$ to minimize $|\setT|$ then pick any one of them.
There are also two simple observations:
\begin{itemize}
\item All rows of the matrix $\mathbf{T}''(\setT^*)$ are still independent.
\item Now all rows of $\mathbf{T}$ are in the row space of $\mathbf{T}''(\setT^*)$.
\end{itemize}

Let $\rho_2$ be the matrix norm of 
$\mathbf{T}(\mathbf{T}''(\setT^*))^\dagger$
induced by the $\ell_1$ vector norm. Let $\epsilon_1$ and $\epsilon_2$ be the same as in \eqref{eq:eps1} and $\eqref{eq:eps2}$, 
and we define the direction vector $\dir\in\Real^{|\setS|}$ as
\begin{align*}
\dir:=&\sgn\left(\mathbf{T}_{\setS}\angh-
\left[\tiny
\begin{array}{c}
(\genulim)_{\setS\cap\VertexG}\\
(\overline{\pflow})_{(\setS-\nGL)\cap\Edge}
\end{array}
\right]\right)\\
&+\sgn\left(\mathbf{T}_{\setS}\angh-
\left[\tiny
\begin{array}{c}
(\genllim)_{\setS\cap\VertexG}\\
(\underline{\pflow})_{(\setS-\nGL)\cap\Edge}
\end{array}
\right]\right)
\end{align*}
where $\sgn$ applies the sign function to each coordinate of the vector.
We then construct 
\begin{align*}
\sload_{(n)}:=\mathbf{T}_{\tiny\VertexL}(\mathbf{T}''(\setT^*))^\dagger
\left[\small
\begin{array}{c}
0\\
\mathbf{T}_{\setS}\angh-\frac{\min\{\epsilon_1,\epsilon_2\}}{2n\nGL\rho_2}\dir\\
\mathbf{T}_{\setR^*}\angh\\
\mathbf{T}_{\setT^*}\angh
\end{array}
\right].
\end{align*}
Since all rows of $\mathbf{T}$ are in the row space of $\mathbf{T}''(\setT^*)$, 
we have
\begin{align*}
(\mathbf{T}''(\setT^*))^\dagger
\left[\tiny
\begin{array}{c}
0\\
\mathbf{T}_{\setS}\angh\\
\mathbf{T}_{\setR^*}\angh\\
\mathbf{T}_{\setT^*}\angh
\end{array}
\right]-\angh
\end{align*}
is perpendicular to the row space of $\mathbf{T}_{\VertexL}$. Therefore,
\begin{align*}
\lim\limits_{n\to\infty}\sload_{(n)}&\to 
\mathbf{T}_{\tiny\VertexL}(\mathbf{T}''(\setT^*))^\dagger
\left[\tiny
\begin{array}{c}
0\\
\mathbf{T}_{\setS}\angh\\
\mathbf{T}_{\setR^*}\angh\\
\mathbf{T}_{\setT^*}\angh
\end{array}
\right]=\mathbf{T}_{\tiny\VertexL}\angh=\sloadh.
\end{align*}
Besides, since
\begin{align*}
\sloadh-\sload_{(n)}=
\mathbf{T}_{\tiny\VertexL}(\mathbf{T}''(\setT^*))^\dagger
\left[\small
\begin{array}{c}
0\\
\frac{\min\{\epsilon_1,\epsilon_2\}}{2n\nGL\rho_2}\dir\\
\zero^{(|\setR^*|+|\setT^*|)\times1}
\end{array}
\right],
\end{align*}
we have
\begin{align*}
\sload_{(n)}&\geq\sloadh-\rho_2|\setS|\frac{\epsilon_2}{2n\nGL\rho_2}\one_{\nL}
\geq\sloadh-\frac{\epsilon_2}{2n}\one_{\nL}>0.
\end{align*}
We then construct the associated $\ang_{(n)}$, $\sgen_{(n)}$ and $\pflow_{(n)}$ as
\begin{align*}
\ang_{(n)}&:=(\mathbf{T}''(\setT^*))^\dagger
\left[\small
\begin{array}{c}
0\\
\mathbf{T}_{\setS}\angh-\frac{\min\{\epsilon_1,\epsilon_2\}}{2n\nGL\rho_2}\dir\\
\mathbf{T}_{\setR^*}\angh\\
\mathbf{T}_{\setT^*}\angh
\end{array}
\right],\\
\sgen_{(n)}&:=\mathbf{T}_{\tiny\VertexG}\ang_{(n)},\\
\pflow_{(n)}&:=\mathbf{T}_{\tiny[\nGL+1,\nGL+\nE]}\ang_{(n)}.
\end{align*}
For $\sgen_{(n)}$, we have
\begin{align*}
|\sgenh-\sgen_{(n)}|&=
\left|
\mathbf{T}_{\tiny\VertexG}(\mathbf{T}''(\setT^*))^\dagger
\left[\small
\begin{array}{c}
0\\
\frac{\min\{\epsilon_1,\epsilon_2\}}{2n\nGL\rho_2}\dir\\
\zero^{(|\setR^*|+|\setT^*|)\times1}
\end{array}
\right]\right|\\
&\leq\rho_2|\setS|\frac{\epsilon_1}{2n\nGL\rho_2}\one_{\nG}\leq\frac{\epsilon_1}{2n}\one_{\nG},
\end{align*}
and consider that all the generators that reach the upper or lower bounds in $\sgenh$ have been moved towards the opposite directions encoded in $\dir$.
All the coordinates in $\sgen_{(n)}$ will then strictly stay within the limits. 
The similar argument also applies to $\pflow_{(n)}$ and implies that all the coordinates in $\pflow_{(n)}$ also strictly stay within the limits.
Thereby, $\sload_{(n)}\in\setsl(\para)$ and there is no binding constraint associated with $(\sgen_{(n)},\sload_{(n)})$.
We have shown in the first case that when no binding constraints arises, there is always an open neighborhood $U(\sload_{(n)})\subseteq\setsl(\para)$. 
We now establish the proof of $\clos(\inte(\setsl(\para)))=\clos(\setsl(\para))$.

Next, we will further show $\setslr(\para)$  is dense in $\setsl(\para)$.
In fact, $\forall\para\in\setparar$, if for some $\sload\in\setsl(\para)$, the optimal solution to \eqref{eq:opf1} has $\geq \nG$ tight inequality constraints, 
then we use $\setS\subseteq[\nGL+\nE]\backslash\VertexL,|\setS|=\nG$ again to denote the indices of any $\nG$ tight inequality constraints.
As those $\nG$ inequality constraints are tight, there must exist $\gamma\in\Gamma$, as defined in \eqref{eq:Gamma}, such that $\gamma^{\T}\para=\mathbf{T}_{\setS}\ang^*$ for the optimal $\ang^*\in\Real^{\nGL}$.
According to \eqref{eq:setupdate}, $\rank(\mathbf{T}_{\setS})$ must be exactly $\nG$. We now have
\begin{subequations}
\begin{eqnarray}
&& \gamma^{\T}\para=\mathbf{T}_{\setS}\ang^*\\
&& -\sload=\mathbf{T}_{\VertexL}\ang^*.
\end{eqnarray}
\label{eq:slplane}
\end{subequations}
For each $\gamma\in\Gamma$, as $\rank(\mathbf{T}_{\setS})=\nG$ but $\rank(\mathbf{T})=\nGL-1$, the set $\{\sload~|~\exists\ang^*,\eqref{eq:slplane}~\text{holds}\}$ is a strict subspace in $\Real^{\nL}$ and thereby nowhere dense in $\setsl$ according to Proposition \ref{setparaissolid} and Lemma \ref{lm:dense}.
As the result, we have
\begin{align*}
\setslr\supseteq\setsl\setminus
\bigcup\limits_{\substack{|\setS|=\nG\\ \setS\subseteq[\nGL+\nE]\backslash\VertexL}}
\bigcup\limits_{\gamma\in\Gamma}
\{\sload~|~\exists\ang^*,\eqref{eq:slplane}~\text{holds for $\gamma$}\}
\end{align*}
must be dense in $\setsl$. 

\section{Proofs of the Lemmas Related to Theorem \ref{thm:range-equiv}}\label{app:proof_wc_reduction}
\begin{IEEEproof}(Lemma \ref{lm:find_f_para1})
We first set $\genllim_i\equiv0$ and $\genulim_i\equiv2$ for all $i\in\VertexG$. Let
\begin{align}\label{eq:construct_sg}
(\sgens)_i=\left\{
\begin{array}{ll}
0, & i\in\setSgen\\
1, & i\not\in\setSgen
\end{array}
\right.
\end{align}
and $\sloads=\frac{\nG-|\setSgen|}{\nL}\one_{\nL}$. 
The construction here guarantees that all $(\sgens)_i$ for $i\in\setSgen$ hit the lower bounds, and other $(\sgens)_i$ are strictly within the bounds.
Then we let
\begin{align*}
\angs=\M^{\dagger}\left[
\begin{array}{c}
\sgens\\
-\sloads\\
0
\end{array}
\right],~\pflows=\B\CM^\T\angs
\end{align*}
where $\M$ is defined in Section \ref{sec:OPF}.
Let
\begin{align*}
&\overline{\pflow}_i=\left\{
\begin{array}{ll}
(\pflows)_i, & {\text{if }} i\in\setSbra~\text{and}~(\pflows)_i\geq0\\
\|\pflows\|_{\infty}+1, & \text{otherwise}
\end{array}
\right.,\\
&\underline{\pflow}_i=\left\{
\begin{array}{ll}
(\pflows)_i, & {\text{if }} i\in\setSbra~\text{and}~(\pflows)_i<0\\
-\|\pflows\|_{\infty}-1, & \text{otherwise}
\end{array}
\right..
\end{align*}
Setting $\paras=[(\genulim)^{\T},(\genllim)^{\T},\overline{\pflow}^{\T},\underline{\pflow}^{\T}]^{\T}$, 
it is easy to check that $(\sgens,\angs)$ is an extreme point of the convex polytope described by \eqref{eq:opf1.b}-\eqref{eq:opf1.e} under $(\paras,\sloads)$ 
since there are exactly $\nGL+\nG$ equality and binding inequality constraints (corresponding to $\setSgen$ and $\setSbra$) in total and they are independent as $\setSgen\perp\setSbra$.
Next, consider the following optimization problem:
\begin{subequations}
\begin{eqnarray}
&\underset{\x}{\text{minimize}}  &\f^{\T}\x
\label{eq:opf_r.a}\\
&\text{ subject to }
&\genllim \leq\x\leq \genulim
\label{eq:opf_r.b}\\
& &\underline{\pflow}\leq\B\CM^{\T}
(\CM\B\CM^{\T})^{\dagger}
\Big[ \begin{array}{c}
\x \\
\nonumber
-\sload
\end{array} \Big]
\leq\overline{\pflow}.
\label{eq:opf_r.c}
\end{eqnarray}
\label{eq:opf_r}
\end{subequations}
Here \eqref{eq:opf1} and \eqref{eq:opf_r} are equivalent to each other in the sense that there is a bijection between their feasible points shown as below.
\begin{align*}
&(\sgen_\fea,\ang_\fea)\to\x_\fea: \x_\fea=\sgen_\fea\\
&\x_\fea\to(\sgen_\fea,\ang_\fea): \sgen_\fea=\x_\fea, \quad \ang_\fea=
\M^{\dagger}\left[
\begin{array}{c}
\sgen_\fea\\
-\sloads\\
0
\end{array}
\right].
\end{align*}
Since $\ang$ is always linear in $\sgen$ for fixed $\sload$, the value of $\sgens$ in \eqref{eq:construct_sg} is also an extreme point of the feasible domain in \eqref{eq:opf_r}.
Therefore there exists $\f'\in\Real^{\nG}$ such that when $\f=\f'$ in \eqref{eq:opf_r}, the optimal solution is uniquely $\x^*=\sgens$.
The equivalence between \eqref{eq:opf1} and \eqref{eq:opf_r} implies that when $\f=\f'$ in \eqref{eq:opf1}, the optimal solution is $(\sgens,\angs)$ and is unique.
Finally, we construct $\fs=\f'+\|\f\|_{\infty}\one$, then the optimal solution remains the same as $(\sgens,\angs)$ and is still unique due to the fact that $\one^\T\sgens\equiv\one^\T\sloads$, but we now have $\fs\geq0$.
\end{IEEEproof}

\begin{IEEEproof}(Lemma \ref{lm:find_f_para2})
We start from $(\fs,\paras,\sloads)$ provided in Lemma \ref{lm:find_f_para1}, and then perturb the parameters in a specific order to derive the desired $(\fss,\parass,W)$.

First, \cite{Zha1990note} shows that the optimal solution set $\OPF(\sloads)$ to \eqref{eq:opf1} for fixed $\sloads$ is both upper hemi-continuous and lower hemi-continuous in $\f$.
Now for the convenience of notation, we use $\OPF^{\f}$ to denote $\OPF(\sloads)$ under the cost vector vector $\f$. For now, $\para$ is chosen to be $\paras$.
Therefore the optimal solution is $(\sgens,\angs)$ and $\OPF^{\fs}=\{\sgens\}$.
As upper hemi-continuity implies that for any neighborhood $U$ of $\sgens$, there is a neighborhood $V$ of $\fs$ such that $\forall \f\in V$, $\OPF^{\f}\subseteq U$.
Consider that \eqref{eq:opf1} is a linear programming problem, so the optimal solution set should contain at least a different extreme point $((\sgen)',\ang')\neq(\sgens,\angs)$ if $\OPF^{\f}\neq\{\sgens\}$.
Here, $(\sgen)'\neq\sgens$ must hold as $(\sgen)'=\sgens$ implies $\ang'=\angs$.
Since a compact convex polytope has only finite extreme points, we can always choose $U$ to be small enough that $(\sgens,\angs)$ is the only extreme point satisfying $\sgens\in U$.
Then there must be a neighborhood $V$ of $\fs$ such that $\forall \f\in V$, $\OPF^{\f}\equiv\{\sgens\}$.
Proposition \ref{fisdense} shows that $\setf$ is dense in $\Real_+^{\nG}$, so there must be some $\fss\in U\cap\setf$ and under $\fss$, 
$\OPF^{\fss}=\{\sgens\}$ and thereby all the binding constraints are the same as the binding constraints under $\fs$, which exactly correspond to $\setSgen$ and $\setSbra$.
In the proof thus far we have taken the parameters in \eqref{eq:opf1} from $(\fs,\paras,\sloads)$ to $(\fss,\paras,\sloads)$.

Next, we are going to perturb $\paras$ to some point in $\setparar(\fss)$. We know that
\begin{itemize}
\item $(\sgens,\angs)$ is the unique solution to \eqref{eq:opf1}.
\item All the constraints and the cost function in \eqref{eq:opf1} are linear and thereby twice continuously differentiable in $(\sgen,\ang)$ and differentiable in $\para$.
\item Since all the binding constraints exactly correspond to $\setSgen$ and $\setSbra$ where $\setSgen\perp\setSbra$, the gradients for all the binding inequalities and equality constraints are independent.
\item We have $|\setSgen|+|\setSbra|=\nG-1$ binding inequality constraints. Together with the fact that $\fss\in\setf$ and thus \eqref{eq:multiplier1} holds, strict complementary slackness must hold.
\end{itemize}
Lemma \ref{Le:existence} shows the set of binding constraints do not change in a small neighborhood $U$ of $\paras$.
Proposition \ref{limitsdense} shows $\setparar(\fss)$ is dense in $\setpara$, so there must be some $\parass\in U\cap\setparar(\fss)$ and under $\parass$, 
all the binding constraints are the same as the binding constraints under $\paras$, which exactly correspond to $\setSgen$ and $\setSbra$.
At this point, the parameters in \eqref{eq:opf1} have been updated to $(\fss,\parass,\sloads)$.

Finally, using the technique similar to the perturbation around $\paras$ above, the set of binding constraints do not change as well when $\sload$ falls within a small neighborhood $U$ of $\sloads$, 
so it is sufficient to show $U\cap\setslr(\parass,\fss)$ contains an open ball $W$. 
First, it is easy to find an open ball $W'$ in $U\cap\setsl(\parass)$ since $\clos(\inte(\setsl(\parass)))=\clos(\setsl(\parass))$ by Proposition \ref{limitsdense} implies that $\sloads$ must be the limit of a sequence of points which are all interior points of $\setsl(\parass)$.
Thus we can always find an interior point of $\setsl(\parass)$ that is strictly within $U$ and take its neighborhood $W'\subseteq U\cap \setsl(\parass)$.
Next, as $\setsl(\parass)\backslash\setslr(\parass,\fss)$ can be covered by the union of finitely many affine hyperplanes, $W'\backslash(\setsl(\parass)\backslash\setslr(\parass,\fss))$ must contain a smaller open ball $W$, which is a subset of $U\cap\setslr(\parass,\fss)$.
\end{IEEEproof}


%\section*{Acknowledgment}

\ifCLASSOPTIONcaptionsoff
  \newpage
\fi



%\bibliographystyle{IEEEtran}
\bibliographystyle{IEEEtran}
\bibliography{references}


%\end{IEEEbiography}


%\newpage



%\vfill


\end{document}